\documentclass[leqno]{article}
\usepackage[utf8]{inputenc}
\usepackage{fontenc}
\usepackage[english]{babel}
\usepackage{amsmath}
\usepackage{amsthm}
\usepackage{systeme}
\usepackage{mathtools, bm}
\usepackage{amssymb, bm}
\usepackage{graphicx}
\usepackage{comment}
\usepackage{xcolor}
\usepackage{textcomp}
\usepackage{framed}
\usepackage{multicol}
\usepackage{geometry}
\definecolor{shadecolor}{gray}{0.9}

\newtheorem{The}{Theorem}[section]

\newtheorem{Ex}{Example}[section]
\newtheorem{Exs}{Examples}[section]
\newtheorem{Prop}{Proposition}[section]
\newtheorem{Lm}{Lemma}[section]
\newtheorem{Rq}{Remark}[section]
\newtheorem{Rqs}{Remarks}[section]
\newtheorem{Cr}{Corollary}[section]
\usepackage{graphicx}
\usepackage{comment}
\usepackage{geometry}
\usepackage[nameinlink,noabbrev]{cleveref}

\newenvironment{Proof}{{\bf Proof:}}
{%
\mbox{}%
\nolinebreak%
\hfill%
\rule{2mm}{2mm}%
\medbreak%
\par%
}

\makeatletter
\def\@seccntformat#1{\@ifundefined{#1@cntformat}%
   {\csname the#1\endcsname\quad}
   {\csname #1@cntformat\endcsname}
}
\makeatother

\title{\bf\large ON THE BIHARMONICITY OF VECTOR FIELDS AND UNIT VECTOR FIELDS }
\author{Mohamed Tahar Kadaoui ABBASSI and
        Souhail DOUA \\ \\
\small{Department of Mathematics, Faculty of Sciences Dhar El Mahraz,} \\
\small{Sidi Mohamed Benabellah University, B.P. 1796-Atlas, Fez, Morocco.}\\
\textit{\footnotesize{E-mail addresses}: mtk{\_}abbassi@Yahoo.fr and souherukun@gmail.com.}}
\date{}

\begin{document}


\maketitle

\begin{abstract}
Let $(M,g)$ be a compact Riemannian manifold. Equipping its tangent bundle $TM$ (resp. unit tangent bundle $T_1M$) by a pseudo-Riemannian $g$-natural metric $G$ (resp. $\tilde{G}$), we study the biharmonicty of vector fields (resp. unit vector fields) as maps $(M,g) \rightarrow (TM,G)$ (resp. $(M,g) \rightarrow (T_1M,\tilde{G})$) as well as critical points of the bienergy functional restricted to the set $\mathfrak{X}(M)$ (resp. $\mathfrak{X}^1(M)$) of vector fields (resp. unit tangent bundles) on $M$. Contrary to the Sasaki metric on $TM$, where the two notions are equivalent to the harmonicity of the vector field and then to its parallelism, we prove that for large classes of $g$-natural metrics on $TM$ the two notions are not equivalent. Furthermore, we give examples of vector fields which are biharmonic as critical points of the bienergy functional restricted to $\mathfrak{X}(M)$, but are not biharmonic maps. We provide equally examples of proper biharmonic vector fields (resp. unit vector fields), i.e. those which are biharmonic without being harmonic.

\end{abstract}

\section{Introduction}


Since its initiation by Eells and Sampson in 1964 (cf.\cite{ES}), the theory of harmonic maps has inspired a lot of mathematicians and physicists to explore interesting problems and to obtain many applications in various domains. It is becoming now among the most widely studied topics in the fields of geometric analysis and differential geometry. A \emph{harmonic map} between two compact Riemannian manifolds $(M,g)$ and $(N,h)$ is a smooth map which is a critical point of the functional energy $E(f)=\frac{1}{2}\int_{M}||df||^{2}v_{g}$. It is characterized by the vanishing of the \emph{tension field} $\tau(f)=tr_{g}\nabla df$. As a generalization of harmonic maps, \emph{biharmonic maps} (cf. \cite{Ee.se}) had been defined as the critical points of the bienergy functional $E_{2}(f)=\frac{1}{2}\int_{M}||\tau(f)||^{2}v_{g}$.\par

In the framework of Riemannian geometry of tangent bundles, particular maps arise naturally, e.g. vector fields, the projection map, tangent maps, etc..., and it is interesting to study their (bi-)harmonicity. Since the (bi-)harmonicity depends on Riemannian structures on the source and target sets of the considered map, the situation would depend evidently, in our case, on the choice of (pseudo-)Riemannian metrics on the tangent bundles involved in the studied maps. Actually, as a (pseudo-)Riemannian manifold, the tangent bundle of a Riemannian manifold had been endowed with various metrics, starting from the \textgravedbl classical\textacutedbl \: Sasaki metric and other lifted metrics from the base manifold, passing through the Cheeger-Gromoll type metrics and the general class of Oproiu metrics and arriving to the wider class of $g$-natural metrics, depending on six independent smooth functions from $\mathbb{R}^{+}$ to $\mathbb{R}$ (cf. \cite{Ab.g}, \cite{Ab.nat}, \cite{Ab.heri} and \cite{11} ). It had been proved that the harmonicity of such particular maps depends on the choice of (pseudo-)Riemannian structures on the involved tangent bundles (cf. \cite{Gilgilala}, \cite{Gil-Med}, \cite{Ish}  for the Sasaki metric and \cite{Ab.unit}, \cite{Ab.har} for $g$-natural metrics).

In the biharmonic setting, M. Markellos and H. Urakawa studied the problem of biharmonicity of vector fields when the tangnet bundle is endowed with the Sasaki metric. Indeed, regarding a tangent vector field on a compact Riemannian manifold $(M,g)$ as a map from $(M,g)$ to $(TM,g^{s})$, where $TM$ is the tangent bundle of $M$ equipped with the Sasaki metric $g^{s}$, they proved that a vector field is biharmonic if and only if it is harmonic, i.e. parallel (cf. \cite{Ur.sect}). They also proved that a vector field is a \emph{biharmonic vector field}, i.e. a critical point of the bienergy functional, restricted to the set $\mathfrak{X}(M)$ of tangent vector fields, if and only if it is parallel.

Such results show that, as in the harmonic case, the  Sasaki metric $g^{s}$ have a kind of \textgravedbl rigidity\textacutedbl regarding the existence of biharmonic vector fields.

So, as it had been done in the harmonic case (cf. \cite{Ab.unit} and \cite{Ab.har}), and to get new examples of biharmonic vector fields in the (non)-compact case other than the harmonic ones, it is worthy to equip the tangent bundles with (pseudo-)Riemannian $g$-natural metrics and to study vector fields as biharmonic maps from $(M,g)$ to $(TM,G)$, as well as \emph{$G$-biharmonic vector fields}, i.e. as critical points of the energy functional restricted  to $\mathfrak{X}(M)$.  \par

On the other hand, C.M. Wood \cite{Wood1} introduced the notion of \emph{harmonic unit vector field} on a compact Riemannian manifold $M $ by restricting the energy functional to the set $\mathfrak{X}^1(M)$ of unit vector fields, $T_1M$ being endowed with the Sasaki metric. Considering arbitrary Riemannian $g$-natural metrics on the unit tangent bundle, corresponding harmonic unit vector fields had been studied in \cite{Ab.unit}.

In the same spirit, M. Markellos and H. Urakawa \cite{Ur.unit} introduced the so-called \emph{biharmonic unit vector fields} as critical points of the bienergy functional restricted to $\mathfrak{X}^1(M)$, the unit tangnet bundle being endowed with the Sasaki metric. It would be interesting to study the question when the unit tangent bundle is equipped with an arbitrary pseudo-Riemannian $g$-natural metric.

In this paper, we will endow the tangent bundle (resp. unit tangent bundle) with an arbitrary pseudo-Riemannian $g$-natural metric $G$ (resp. $\tilde{G}$), and we will study the biharmonicity of a vector (resp. unit vector) fields, as maps from $(M,g)$ to $(TM,G)$ (resp. $(T_1M,\tilde{G})$) and as critical points of the bienergy functional restricted to $\mathfrak{X}(M)$ (resp. $\mathfrak{X}^1(M)$).  More precisely, we will be interested to answer to the following questions:
\begin{enumerate}
  \item Are there parallel vector fields on $M$ which are not biharmonic maps from $(M,g)$ to $(TM,G)$?
  \item Are there vector fields on $M$ which are \emph{proper biharmonic}, i.e. biharmonic without being harmonic?
  \item Are there vector fields on $M$ which are \emph{$G$-biharmonic} (i.e. a critical point of the bienergy functional restricted to $\mathfrak{X}(M)$), but not biharmonic?
  \item Give a tensorial characterization of \emph{$\tilde{G}$-biharmonic unit vector fields}, i.e. critical point of the bienergy functional restricted to $\mathfrak{X}^1(M)$.
  \item Are there unit vector fields on $M$ which are \emph{proper} $\tilde{G}$-biharmonic, i.e. $\tilde{G}$-biharmonic but not harmonic?
\end{enumerate}

The paper is organised as follows. In section 2, we give a brief exposition of some of the basic notions about the geometry of tangent bundle and biharmonic maps. As mentioned above, when equipping the tangent bundle of a compact Riemannian manifold by the Sasaki metric, the only vector fields which are biharmonic, as maps, are the parallel ones. One can wonder if this is still true if we equip the tangent bundle with an arbitrary pseudo-Riemannian $g$-natural metric. As it is very hard to deal with the general case, we will restrict ourselves to the study of parallel vector fields and we devote section 3 to compute the bitension field of a parallel vector field and to give necessary and sufficient conditions for a parallel vector field to be biharmonic (Theorem \ref{gbitension_field} for the general case and Theorem \ref{biharm-KKtype} and Corollary \ref{Biharm-KKtype2} for the special case of Kaluza-Klein type metrics). In particular, we answer to the two first questions by giving examples of non biharmonic parallel vector fields (Remark \ref{Rq-biharm}) and of parallel proper biharmonic vector fields (Examples \ref{Exs-prop-biharm}).

In section 4, we (will) derive the condition for a parallel vector field to be a $G$-biharmonic vector field (Theorem \ref{G-0}). As a consequence, we give a subclass of pseudo-Riemannian $g$-natural metrics on $TM$ for which the $G$-biharmonicity of parallel vector fields is equivalent to their biharmonicity (Example \ref{Ex-equiv}), generalizing a same result for the Sasaki metric. As an answer to the third question, we construct a 2-parameter family of pseudo-Riemannian $g$-natural metrics on $TM$ for which every parallel vector field is $G$-biharmonic but not biharmonic (Proposition \ref{G-biharm-not-biharm} and Example \ref{Ex-G-biharm-not-biharm}). \par

Section 5 is devoted to the characterization of $\tilde{G}$-biharmonic unit vector fields, i.e. critical points of the bienergy functional restricted to $\mathfrak{X}^1(M)$, with applications to special classes of vector fields e.g. Killing vector fields and Reeb vector fields. We give examples of proper $G$-biharmonic unit vector fields on hyperbolic spaces $\mathbb{H}^n$ (Example \ref{Ex-hyperbolic}), on the solvable Lie group $Sol_3$ (Example \ref{Ex-Solvable}) and on the special unitary group $SU(2)$ with a left invariant Riemannian metric (Example \ref{Ex-Sp-unitary}).

Finally, we would like to thank E. Loubeau for his valuable comments and suggestions on a first version of this paper.


\section{Preliminaries}


\subsection{Harmonic and biharmonic maps}


Let $(M,g)$ be an $m$-dimensional compact Riemannian manifold and $(N,h)$ be an $n$-dimensional pseudo-Riemannian manifold. If $f: (M,g) \rightarrow (N,h)$ is a smooth map , then the  {\em energy} of $f$ is defined as the integral
$$ E(f) := \int _M e(f) dv_g $$
where $e(f)= \frac{1}{2} ||f_*||^2= \frac{1}{2} {\rm tr}_g f^* h$ is  the so-called {\em energy density} of $f$. With respect to a local orthonormal frame $\{e_1,..,e_n \}$ on $M$, it is possible to express the energy density as $e(f)= \frac{1}{2} \sum _{i=1} ^n h(f_*e_i,f_*e_i)$. Critical points of the energy functional on $C^{\infty} (M,N)$ are known as {\em harmonic maps}. They have been characterized in \cite{ES} as maps having vanishing {\em tension field} $\tau (f)={\rm tr} \nabla df$. When $(M,g)$ is a general Riemannian manifold (including the non-compact case), a map $f: (M,g) \rightarrow (N,h)$ is said to be \emph{harmonic} if $\tau (f)=0$. For further details about harmonic maps, we can refer to \cite{report}.

In \cite{Ee.se}, J. Eelles and L. Lemaire extended  the notion of harmonic maps to \emph{biharmonic maps}, which are, by definition, critical points of the \emph{bienergy funtional}:
$$E_{2}(f)=\frac{1}{2}\int_{M} ||\tau(f)||^{2} dv_{g}.$$
G. Jiang (\cite{Ji}) derived the associated Euler-Lagrange equation of $E_{2}$ as
\begin{equation}\label{bitension_field}
\widehat{\Delta}\tau(f)-\underset{i=1}{\overset{m}{\sum}}R^{N}\big(\tau(f),f_{*}e_{i}\big)f_{*}e_{i}=0,
\end{equation}
where $R^{N}$ is the Riemannian curvature tensor of $(N,h)$ defined by: $R^{N}(X,Y)Z=\nabla^{N}_{X}\nabla^{N}_{Y}Z-\nabla^{N}_{Y}\nabla^{N}_{X}Z-\nabla^{N}_{[X,Y]}Z$ for $X$, $Y$, $Z$ $\in \mathfrak{X}(N)$ and $\widehat{\Delta}$ is the rough laplacien defined by: $\widehat{\Delta}=-\underset{i=1}{\overset{m}{\sum}}\big(\widehat{\nabla}_{e_{i}}\widehat{\nabla}_{e_{i}}-\widehat{\nabla}_{\nabla_{e_{i}}e_{i}}\big)$, where $\widehat{\nabla}$ is the induced connexion on $f^{-1}TN$.\par

The quantity $\tau_{2}(f):=\widehat{\Delta}\tau(f)-\underset{i=1}{\overset{m}{\sum}}R^{N}\big(\tau(f),f_{*}e_{i}\big)f_{*}e_{i}$ is called the \emph{bitension field of $f$}. As for the case of harmonic maps, when $(M,g)$ is a general Riemannian manifold (including the non-compact case), a map $f: (M,g) \rightarrow (N,h)$ is said to be \emph{biharmonic} if $\tau_2 (f)=0$.


\subsection{Geometry of tangent bundle}


\subsubsection{Basic formulas on tangent bundles}


Let $\nabla$ be the Levi-Civita connection of $g$. Then the tangent space of $TM$ at any point $(x,u)\in TM$ splits into the horizontal and vertical subspaces with respect to $\nabla$:
$$(TM)_{(x,u)}=H_{(x,u)}\oplus V_{(x,u)}.$$    \par

If $(x,u)\in TM$ is given then, for any vector $X\in M_x$, there exists a unique vector $X^h  \in H_{(x,u)}$ such that $p_* X^h =X$. We call $X^h$ the \emph{horizontal lift} of $X$ to the point $(x,u)\in TM$. The \emph{vertical lift} of a vector $X\in M_x$ to $(x,u)\in TM$ is a vector $X^v  \in V_{(x,u)}$ such that $X^v (df) =Xf$, for all functions $f$ on $M$. Here we consider $1$-forms $df$ on $M$ as functions on $TM$ (i.e. $(df)(x,u)=uf$). Note that the map $X \to X^h$ is an isomorphism between the vector spaces $M_x$ and $H_{(x,u)}$. Similarly, the map $X \to X^v$ is an isomorphism between the vector spaces $M_x$ and $V_{(x,u)}$. Obviously, each tangent vector $\tilde Z \in (TM)_{(x,u)}$ can be    written in the form $\tilde Z =X^h + Y^v$, where $X,Y \in M_x$ are uniquely determined vectors.\par

Let $T$ be a tensor field of type $(1,s)$ on $M$. If $X_1$, $X_2$,..., $X_{s-1}$ $\in M_x$, then $h\{T(X_1,...,u,...,X_{s-1})\}$ (resp.$v\{T(X_1,...,u,...,X_{s-1})\}$) is a horizontal (resp.vertical) vector at $(x,u)$ which is introduced by the formula
\begin{equation*}
h\{T(X_1,..., u,...,X_{s-1})\}=\sum u^\lambda (T(X_1,..., \left(\frac {\partial} {\partial x^\lambda}\right)_x,...,X_{s-1}))^h
\end{equation*}
\begin{equation*}
(\textup{resp.}v\{T(X_1,..., u,...,X_{s-1})\}=\sum u^\lambda(T(X_1,...,\left(\frac {\partial} {\partial x^\lambda}\right)_x,...,X_{s-1}))^v).
\end{equation*}
In particular, if $T$ is the identity tensor of type $(1,1)$, then we obtain the geodesic flow vector field at $(x,u)$, $\xi_{(x,u)}=\sum u^\lambda (\frac {\partial} {\partial
x^\lambda})^h_{(x,u)}$, and the canonical vertical vector at $(x,u)$, $\mathcal{U}_{(x,u)}=\sum u^\lambda (\frac {\partial}{\partial x^\lambda})^v_{(x,u)}$.\\ Moreover
$h\{T(X_1,..., u,...,u,...,X_{s-t})\}$ and $v\{T(X_1,..., u,...,u,...,X_{s-t})\}$ are introduced by similar way.\\ Also we make the conventions
\begin{equation*}
h\{T(X_1,...,X_{s})\}=(T(X_1,...,X_{s}))^h \; \mbox{and} \;
v\{T(X_1,...,X_{s})\}=(T(X_1,...,X_{s}))^v.
\end{equation*}
Thus $h\{X\}=X^h$ and $v\{X\}=X^v$, for each vector field $X$ on $M$.\par

From the preceding quantities, one can define vector fields on $TU$ in the following way: If $u= \sum_i u^i \left(\frac {\partial}{\partial x^i}\right)_x$ is a fixed point in $TM$ and $X_1,...,X_{s-1}$ are vector fields on $U$, then we denote by
\begin{equation*}
h\{T(X_1,..., u,...,X_{s-1})\} \; \mbox{(resp.} \; v\{T(X_1,..., u,...,X_{s-1})\} \mbox{)}
\end{equation*}
the horizontal (resp. vertical) vector field on $TU$ defined by
\begin{equation*}
h\{T(X_1,...,u,...,X_{s-1})\} =\sum_{\lambda} u^{\lambda} [T(X_1,..., \frac {\partial} {\partial x^{\lambda}},...,X_{s-1})]^h
\end{equation*}
\begin{equation*}
(\mbox{resp.} \; v\{T(X_1,...,u,...,X_{s-1})\} =\sum_{\lambda} u^{\lambda} [T(X_1,..., \frac {\partial} {\partial x^{\lambda}},...,X_{s-1})]^v.
\end{equation*}
Moreover, for vector fields $X_1,...,X_{s-1}$ on $U$, the vector fields $h\{T(X_1,..., u,...,$ $u,...,X_{s-t})\}$ and $v\{T(X_1,..., u,...,u,...,X_{s-t})\}$, on $TU$, are introduced by similar way.


\subsubsection{$g$-natural metrics}


All $g$-natural metrics on the tangent bundle of a Riemannian manifold $(M,g)$ are completely determined as follows:

\begin{Prop}\label{$g$-nat} \cite{Ab.nat}
Let $(M,g)$ be a Riemannian manifold and $G$ be a $g$-natural metric on $TM$. Then there are functions $\alpha_i$, $\beta_i:\mathbb{R}^+ \rightarrow \mathbb{R}$, $i=1,2,3$, such     that for every $u$, $X$, $Y\in M_x$, we have
    \arraycolsep1.5pt
    \begin{equation}\label{g_exp}
     \left\lbrace
    \begin{array}{rcl}
    G_{(x,u)}(X^h,Y^h)& = & (\alpha_1+ \alpha_3)(\rho) g_x(X,Y)
                        + (\beta_1+ \beta_3)(\rho)g_x(X,u)g_x(Y,u),\\
    G_{(x,u)}(X^h,Y^v)& = & \alpha_2 (\rho) g_x(X,Y)
                           +  \beta_2 (\rho) g_x(X,u)g_x(Y,u), \\
    G_{(x,u)}(X^v,Y^h)& = & \alpha_2 (\rho) g_x(X,Y)
                           +  \beta_2 (\rho) g_x(X,u)g_x(Y,u), \\
    G_{(x,u)}(X^v,Y^v)& = & \alpha_1 (\rho) g_x(X,Y)
                           + \beta_1 (\rho) g_x(X,u)g_x(Y,u),
    \end{array}
        \arraycolsep5pt \right.
    \end{equation}
where $\rho =g_x(u,u)$.
\end{Prop}

\begin{Rq}
From now on, we shall use the following notations:
    \begin{itemize}
      \item $\phi_i(t) =\alpha_i(t) +t \beta_i(t)$,
      \item $\alpha(t) = \alpha_1(t) (\alpha_1+\alpha_3)(t) - \alpha_2 ^2$(t),
      \item $\phi(t) = \phi_1(t) (\phi_1 +\phi_3)(t) -\phi_2^2(t)$,
    \end{itemize}
for all $t \in \mathbb{R}^+$.
    \end{Rq}

Pseudo-Riemannian $g$-natural metrics are characterized as follows:

\begin{Prop}\label{riem-nat} \cite{Ab.nat}
A $g$-natural metric $G$ on the tangent bundle of a Riemannian manifold $(M,g)$, defined by the functions of Proposition \ref{$g$-nat}, is
\begin{itemize}
	\item non-degenerate if and only if
	$$\alpha(t) \neq 0, \qquad  \phi(t) \neq 0 \qquad \text{for all} \;\, t \in \mathbb{R}^+ ; $$
	\item Riemannian if and only if
	$$\alpha_1(t) > 0,  \qquad \phi_1(t) > 0, \qquad \alpha(t) > 0, \qquad  \phi(t)>0 \qquad \text{for all} \;\, t \in \mathbb{R}^+ . $$
\end{itemize}
for all $t \in \mathbb{R}^+$.\par
\end{Prop}

The wide class of $g$-natural metrics includes several well known metrics (Riemannian and not) on $TM$. In particular:
\begin{itemize}
	\item the {\em Sasaki metric} $g_S$ is obtained for $\alpha _1 =1$ and $\alpha _2 = \alpha _3 = \beta _1 =\beta _2 = \beta _3 =0$.
	\item {\em Kaluza--Klein metrics}, as commonly defined on principal bundles \cite{Wood1} {(see also \cite{Loubeau})}, are obtained for
	$\alpha _2 = \beta _2 = \beta _1 +\beta _3 = 0$.
	\item {\em Metrics of Kaluza--Klein type} are defined by the geometric condition of orthogonality between horizontal and vertical distributions. Thus, a  $g$-natural metric $G$ is of Kaluza-Klein type if $\alpha _2=\beta _2 =0$.
\end{itemize}

\begin{Rqs}
\begin{enumerate}
  \item In the sequel, when we consider an arbitrary $g$-natural metric $G$ on $TM$, we implicitly suppose that it is defined by the functions $\alpha_i$, $\beta_i:\mathbb{R}^+ \rightarrow \mathbb{R}$, $i=1,2,3$, given in Proposition \ref{$g$-nat}.
  \item Unless otherwise stated, all real functions $\alpha_i$, $\beta_i$, $\phi_i$, $\alpha$ and $\phi$ and their derivatives are evaluated at $\rho:=g_x(u,u)$.
\end{enumerate}
\end{Rqs}


\subsubsection{The Levi-Civita connexion of pseudo-Riemannian $g$-natural metrics}


\begin{Prop}[]\label{lev-civ-con}
Let $(M,g)$ be a Riemannian manifold, $\nabla$ its Levi-Civita connection and $R$ its curvature tensor. Let $G$ be a pseudo-Riemannian $g$-natural metric on $TM$. Then the Levi-Civita connection $\overline{\nabla}$ of $(TM,G)$ is characterized by
     \arraycolsep1.5pt
    \begin{equation}\label{cnx_exp}
    \left\lbrace
    \begin{array}{rcl}
    (\overline{\nabla}_{X^h}Y^h)_{(x,u)} & = &  (\nabla_X Y)^h _{(x,u)} +h\{A(u;X_x,Y_x)\} + v\{B(u;X_x,Y_x)\}, \\
     (\overline{\nabla}_{X^h}Y^v)_{(x,u)} & = &  (\nabla_X Y)^v _{(x,u)} +h\{C(u;X_x,Y_x)\} + v\{D(u;X_x,Y_x)\}, \\
     (\overline{\nabla}_{X^v}Y^h)_{(x,u)} & = &  h\{C(u;Y_x,X_x)\} + v\{D(u;Y_x,X_x)\}, \\
     (\overline{\nabla}_{X^v}Y^v)_{(x,u)} & = &  h\{E(u;X_x,Y_x)\} + v\{F(u;X_x,Y_x)\},
    \end{array}  \right.
    \end{equation}
for all vector fields $X$, $Y$ on $M$ and $(x,u) \in TM$, where $A$, $B$, $C$, $D$, $E$ and $F$ are $F$-tensor fields of type $(1,2)$ on $M$ (cf. the definition and some properties on $F$-tensor fields in Appendix A) defined, for all $u$, $X$, $Y \in M_x$, $x \in M$, by:
    $$
    \begin{array}{lcl}
    A(u;X,Y)  & = & A_1 [R(X,u)Y +R(Y,u)X]+A_2 [g_x(Y,u)X + g_x(X,u)Y] \vphantom{\displaystyle\frac{a}{a}} \\
    & &+A_3 g_x(R(X,u)Y,u)u + A_4  g_x(X,Y)u + A_5 g_x(X,u)g_x(Y,u)u,
    \vphantom{\displaystyle\frac{a}{a}}
    \end{array}
    \arraycolsep5pt$$
    where
    \begin{equation}\label{Ai}
    \begin{array}{lcl}
    \vphantom{\displaystyle\frac{A}{A}} A_1 &=& -\frac{\alpha_1 \alpha_2}{2\alpha} , \\
    \vphantom{\displaystyle\frac{A}{A}} A_2 &=& \frac{\alpha_2 (\beta_1 +\beta_3)}{2\alpha}, \\
    \vphantom{\displaystyle\frac{A}{A}} A_3 &=& \frac{ \alpha_2\{ \alpha_1[\phi_1 (\beta_1 +\beta_3) -\phi_2 \beta_2] +\alpha_2(\beta_1 \alpha_2
    - \beta_2 \alpha_1)\}}{\alpha \phi} , \\
    \vphantom{\displaystyle\frac{A}{A}} A_4 &=& \frac{ \phi_2 (\alpha_1 + \alpha_3)^\prime }{\phi} , \\
    \vphantom{\displaystyle\frac{A}{A}} A_5 &=& \frac{\alpha \phi_2(\beta_1 +\beta_3)^\prime + (\beta_1 +\beta_3)\{ \alpha_2 [\phi_2 \beta_2 - \phi_1 (\beta_1 +\beta_3)] + (\alpha_1 +\alpha_3)(\alpha_1 \beta_2 - \alpha_2 \beta_1)\}}{\alpha \phi}, \vphantom{\displaystyle\frac{A}{A}}
    \end{array}
    \arraycolsep5pt
    \end{equation}

    $$\arraycolsep1.5pt
    \begin{array}{lcl}
    B(u;X,Y) & = &  B_1 R(X,u)Y +B_2 R(X,Y)u +B_3 [g_x(Y,u)X +g_x(X,u)Y] \vphantom{\displaystyle\frac{a}{a}} \\
    & & + B_4 g_x(R(X,u)Y,u)u + B_5 g_x(X,Y)u +B_6 g_x(X,u)g_x(Y,u)u , \vphantom{\displaystyle\frac{a}{a}}
    \end{array}
    \arraycolsep5pt$$
    where
    \begin{equation}\label{Bi}
\arraycolsep1.5pt
\begin{array}{lcl}
\vphantom{\displaystyle\frac{A}{A}} B_1 &=&  \frac{\alpha_2 ^2}{\alpha}, \\
\vphantom{\displaystyle\frac{A}{A}} B_3 &=& - \frac{(\alpha_1 +\alpha_3) (\beta_1 +\beta_3)} {2\alpha}, \\
\vphantom{\displaystyle\frac{A}{A}} B_4 &=&  \frac{\alpha_2 \{ \alpha_2[\phi_2 \beta_2- \phi_1 (\beta_1
               +\beta_3)] +(\alpha_1 +\alpha_3)(\beta_2
               \alpha_1  - \beta_1 \alpha_2)\}}{\alpha\phi}, \\
\vphantom{\displaystyle\frac{A}{A}} B_5 &=& - \frac{(\phi_1 +\phi_3)(\alpha_1+ \alpha_3)^\prime}{\phi}, \\
\vphantom{\displaystyle\frac{A}{A}} B_6 &=& \frac{-\alpha (\phi_1 +\phi_3)(\beta_1 +\beta_3)^\prime +(\beta_1 +\beta_3) \{ (\alpha_1 +\alpha_3)[(\phi_1 +\phi_3)
                \beta_1- \phi_2 \beta_2] 
                + \alpha_2[\alpha_2 (\beta_1 +\beta_3)
                -(\alpha_1 +\alpha_3) \beta_2 ] \}}{\alpha\phi} ,
\end{array}
\arraycolsep5pt
\end{equation}

    $$\arraycolsep1.5pt
    \begin{array}{lcl}
    C(u;X,Y) & = &  C_1 R(Y,u)X +C_2 g_x(X,u)Y +C_3 g_x(Y,u)X +C_4 g_x(R(X,u)Y,u)u+ \vphantom{\displaystyle\frac{A}{A}} \vphantom{\displaystyle\frac{A}{A}} \\
                        & & +C_5 g_x(X,Y)u +C_6 g_x(X,u)g_x(Y,u) u ,
    \end{array}
    \arraycolsep5pt$$

    where
    \begin{equation}\label{Ci}
    \arraycolsep1.5pt
    \begin{array}{lcl}
    \vphantom{\displaystyle\frac{A}{A}} C_1 &=& -\frac{\alpha_1^2}{2 \alpha}, \\
    \vphantom{\displaystyle\frac{A}{A}} C_2 &=& -\frac{\alpha_1 (\beta_1 +\beta_3)}{2 \alpha}, \\
    \vphantom{\displaystyle\frac{A}{A}} C_3 &=&  \frac{\alpha_1 (\alpha_1 +\alpha_3)^\prime -\alpha_2(\alpha_2 ^\prime -\frac{\beta_2}{2})}{\alpha}, \\
    \vphantom{\displaystyle\frac{A}{A}} C_4 &=&  \frac{ \alpha_1 \{ \alpha_2 (\alpha_2 \beta_1 -\alpha_1 \beta_2) +\alpha_1 [\phi_1 (\beta_1 +\beta_3) -\phi_2 \beta_2]\}}{2\alpha \phi}, \\
    \vphantom{\displaystyle\frac{A}{A}} C_5 &=&  \frac{\phi_1 (\beta_1 +\beta_3) +\phi_2(2 \alpha_2 ^\prime -\beta_2)}{2\phi} , \\
    \vphantom{\displaystyle\frac{A}{A}} C_6 &=&  \frac{\alpha \phi_1  (\beta_1 +\beta_3)^\prime + \{ \alpha_2(\alpha_1 \beta_2 -\alpha_2 \beta_1) +\alpha_1[\phi_2 \beta_2 -(\beta_1 +\beta_3) \phi_1] \} [(\alpha_1 +\alpha_3)^\prime +\frac{\beta_1 +\beta_3}2] }{\alpha \phi} \\
     & & + \frac{\{\alpha_2 [\beta_1 (\phi_1 +\phi_3) -\beta_2 \phi_2] -\alpha_1 [\beta_2 (\alpha_1 +\alpha_3)  -\alpha_2(\beta_1 +\beta_3)]\}(\alpha_2 ^\prime
     -\frac{\beta_2}{2})}{\alpha \phi},
    \end{array}
    \arraycolsep5pt
    \end{equation}

    $$\arraycolsep1.5pt
    \begin{array}{lcl}
    D(u;X,Y) & = & D_1 R(Y,u)X +D_2 g_x(X,u)Y + D_3 g_x(Y,u)X +D_4 g_x(R(X,u)Y,u)u \vphantom{\displaystyle\frac{A}{A}} \\
             &  &+ D_5 g_x(X,Y)u +D_6g_x(X,u)g_x(Y,u) u ,
    \end{array}
    \arraycolsep5pt$$

    where
    \begin{equation}\label{Di}
    \arraycolsep1.5pt
    \begin{array}{lcl}
    \vphantom{\displaystyle\frac{A}{A}} D_1 &=& \frac{\alpha_1 \alpha_2}{2 \alpha}, \\
    \vphantom{\displaystyle\frac{A}{A}} D_2 &=&  \frac{\alpha_2 (\beta_1 +\beta_3)}{2 \alpha}, \\
    \vphantom{\displaystyle\frac{A}{A}} D_3 &=& \frac{-\alpha_2  (\alpha_1 +\alpha_3)^\prime +(\alpha_1 +\alpha_3) (\alpha_2 ^\prime -\frac{\beta_2}{2})}{\alpha}, \\
    \vphantom{\displaystyle\frac{A}{A}} D_4 &=&  \frac{\alpha_1 \{ (\alpha_1 +\alpha_3)(\alpha_1 \beta_2 -\alpha_2 \beta_1) +\alpha_2 [\phi_2 \beta_2
                  -\phi_1 (\beta_1+\beta_3)]\}}{2\alpha \phi}, \\
          \vphantom{\displaystyle\frac{A}{A}} D_5 &=& -\frac{\phi_2 (\beta_1 +\beta_3) +(\phi_1 +\phi_3)(2\alpha_2 ^\prime -\beta_2)}{2\phi} , \\
    \vphantom{\displaystyle\frac{A}{A}} D_6 &=& \frac{-\alpha \phi_2 (\beta_1 +\beta_3)^\prime  +\{ (\alpha_1 +\alpha_3)(\alpha_2 \beta_1
                   -\alpha_1 \beta_2) +\alpha_2[\phi_1(\beta_1 +\beta_3) -\phi_2 \beta_2]\} [(\alpha_1 +\alpha_3)^\prime +\frac{\beta_1 +\beta_3}2]}{\alpha \phi} \\
     & & + \frac{\{ (\alpha_1 +\alpha_3)[\beta_2 \phi_2 -\beta_1 (\phi_1 +\phi_3)]+\alpha_2 [\beta_2 (\alpha_1 +\alpha_3) -\alpha_2(\beta_1 +\beta_3)]\}(\alpha_2 ^\prime
                   -\frac{\beta_2}{2})}{\alpha \phi},
    \end{array}
    \arraycolsep5pt
    \end{equation}

    $$\arraycolsep1.5pt
    \begin{array}{lcl}
    E(u;X,Y) & = &  E_1 [g_x(Y,u)X +g_x(X,u)Y] + E_2 g_x(X,Y)u +E_3 g_x(X,u)g_x(Y,u) u,
    \end{array}
    \arraycolsep5pt$$

    where
    \begin{equation}\label{Ei}
    \arraycolsep1.5pt
    \begin{array}{lcl}
    \vphantom{\displaystyle\frac{A}{A}} E_1 &=&  \frac{\alpha_1 (\alpha_2^\prime +\frac{\beta_2}{2}) -\alpha_2 \alpha_1^\prime}{\alpha} , \\
    \vphantom{\displaystyle\frac{A}{A}} E_2 &=&  \frac{\phi_1 \beta_2 -\phi_2(\beta_1 -\alpha_1 ^\prime)}{\phi} , \\
    \vphantom{\displaystyle\frac{A}{A}} E_3 &=&  \frac{\alpha (2\phi_1 \beta_2^\prime -\phi_2 \beta_1^\prime) +2\alpha_1^\prime \{ \alpha_1[\alpha_2 (\beta_1 +\beta_3) -
                   \beta_2 (\alpha_1 +\alpha_3)] +\alpha_2[\beta_1(\phi_1 +\phi_3) -\beta_2 \phi_2]\} }{\alpha \phi} \\
                    & & + \frac{(2\alpha_2 ^\prime +\beta_2) \{ \alpha_1[\phi_2 \beta_2 -\phi_1(\beta_1 +\beta_3)] +\alpha_2 (\alpha_1 \beta_2 -\alpha_2 \beta_1)\}}{\alpha \phi},
    \end{array}
    \arraycolsep5pt
    \end{equation}

    $$\arraycolsep1.5pt
    \begin{array}{lcl}
    F(u;X,Y) & = &  F_1 [g_x(Y,u)X +g_x(X,u)Y] +F_2 g_x(X,Y)u +F_3 g_x(X,u)g_x(Y,u) u,
    \end{array}
    \arraycolsep5pt$$

    where
    \begin{equation}\label{Fi}
    \arraycolsep1.5pt
    \begin{array}{lcl}
    \vphantom{\displaystyle\frac{A}{A}} F_1 &=& \frac{-\alpha_2 (\alpha_2^\prime +\frac{\beta_2}{2}) +(\alpha_1 +\alpha_3) \alpha_1^\prime}{\alpha}, \\
    \vphantom{\displaystyle\frac{A}{A}} F_2 &=& \frac{ (\phi_1 +\phi_3)(\beta_1 -\alpha_1 ^\prime) -\phi_2 \beta_2}{\phi}, \\
    \vphantom{\displaystyle\frac{A}{A}} F_3 &=& \frac{\alpha [(\phi_1 +\phi_3) \beta_1^\prime -2\phi_2 \beta_2^\prime] +2\alpha_1^\prime \{ \alpha_2
                   [\beta_2 (\alpha_1 +\alpha_3) - \alpha_2 (\beta_1 +\beta_3)] + (\alpha_1 +\alpha_3) [\beta_2 \phi_2 -\beta_1(\phi_1+\phi_3)]\}}{\alpha \phi} \\
                    & & + \frac{(2\alpha_2 ^\prime +\beta_2) \{ \alpha_2[\phi_1(\beta_1 +\beta_3) -\phi_2 \beta_2] +(\alpha_1 +\alpha_3)(\alpha_2\beta_1
                   -\alpha_1 \beta_2)\}}{\alpha \phi}.
    \end{array}
    \arraycolsep5pt
    \end{equation}
\end{Prop}


\subsubsection{Riemaniann curvature of $g$-natural metrics}


\begin{Prop}\cite{Ab.heri}\label{curvature_exp}
Let $(M,g)$ be a Riemannian manifold and $G$ be a pseudo-Riemannian $g$-natural metric on $TM$. Denote by $\nabla$ and $R$ the Levi-Civita connection and the Rieamannian curvature tensor of $(M,g)$, respectively. Then we have:
     \begin{equation*}
       \begin{array}{rcl}
    \overset{-}{R}\big(X^{h},Y^{h}\big)Z^{h}&=& \Big[R(X,Y)Z\Big]^{h}+h\Bigg\{(\nabla_{X}A_{u})(Y,Z)-(\nabla_{Y}A_{u})(X,Z)  \\

                     &&+A\Big(u;X,A(u;Y,Z)\Big)-A\Big(u;Y,A(u;X,Z)\Big)+C\Big(u;X,B(u;Y,Z)\Big)\\
                     &&-C\Big(u;Y,B(u;X,Z)\Big)+C\Big(u;Z,R(X,Y)u\Big) \Bigg\}\\
                     && +v\Bigg\{(\nabla_{X}B_{u})(Y,Z)-(\nabla_{Y}A_{u})(X,Z)   \\
                     &&+B\Big(u;X,A(u;Y,Z)\Big)-B\Big(u;Y,A(u;X,Z)\Big)+D\Big(u;X,B(u;Y,Z)\Big)\\

                    &&-D\Big(u;Y,B(u;X,Z)\Big) +D\Big(u;Z,R(X,Y)u\Big)    \Bigg\},\\
    \end{array}
    \end{equation*}
    \begin{equation*}
       \begin{array}{rcl}
    \overset{-}{R}\big(X^{h},Y^{v}\big)Z^{h}&=& h\Bigg\{(\nabla_{X}C_{u})(Z,Y)+A\Big(u;X,C(u;Z,Y)\Big)+C\Big(u;X,D(u;Z,Y)\Big)\\
                   &&-C\Big(u;A(u;X,Z),Y\Big)-E\Big(u;Y,B(u;X,Z)\Big)
                    -d_u\big(A_{(X,Z)}\big)(Y)\Bigg\}  \\
                    &&+v\Bigg\{(\nabla_{X}D_{u})(Z,Y)+B\Big(u;X,C(u;Z,Y)\Big)+D\Big(u;X,D(u;Z,Y)\Big)\\
                    &&-D\Big(u;A(u;X,Z),Y\Big)-F\Big(u;Y,B(u;X,Z)\Big)-d_u\Big(B_{(X,Z)}\Big)(Y)\Bigg\}, \\

    \end{array}
    \end{equation*}
where $A_{(X,Z)}$ and $B_{(X,Z)}$ are the mappings $M_x \rightarrow M_x$ defined by \eqref{map-ass} and $d_u\big(A_{(X,Z)}\big)$ and $d_u\big(A_{(X,Z)}\big)$ are their derivatives at $u$, respectively (cf. Appendix A for more details).
\end{Prop}


\section{The biharmonicity of parallel vector fields}


Let $(M,g)$ be a Reimannian manifold, and $TM$ its tangent bundle which is equipped with a pseudo-Riemannian $g$-natural metric $G$. Let $V$ be a vector field on $M$, then $V$ can be regarded as a map $V:(M,g)\longrightarrow (TM,G)$. The bienergy $E_{2}(V)$ of $V$ is the bienergy associated to the map $V:(M,g)\longrightarrow (TM,G)$.\par

The tension field $\tau(V)$ and bitension field $\tau_2(V)$ of $V$ can be decomposed as follows: For all $x \in M$, we denote by $\tau_h(V)(x)$, $\tau_v(V)(x)$, $\tau_{2h}(V)(x)$ and $\tau_{2v}(V)(x)$ the vectors of $M_x$, such that
$$(\tau(V))(x)=(\tau_{h}(V)(x))^h +(\tau_{v}(V)(x))^v, \qquad (\tau_{2}(V))(x)=(\tau_{2h}(V)(x))^h +(\tau_{2v}(V)(x))^v,$$
obtaining four vector fields $\tau_h(V)$, $\tau_v(V)$, $\tau_{2h}(V)$ and $\tau_{2v}(V)$ on $M$. So, as sections in $\Gamma(V^{-1} TTM)$, the tension and bitension fields can be written as
$$(\tau(V))=(\tau_{h}(V))^h \circ V +(\tau_{v}(V))^v \circ V, \qquad (\tau_{2}(V))=(\tau_{2h}(V))^h \circ V +(\tau_{2v}(V))^v \circ V.$$


\subsection{The bitension field of a parallel vector field}


In \cite{Ab.har}, the authors had calculated the tension field of a parallel vector field $V$. More precisely, they found that
     \begin{equation}\label{tenspar1}
     \tau _h (V)=\big(2A_2 +mA_4  +\rho A_5 \big)(\rho)V ,\qquad      \tau _v (V)= \big(2B_3 +mB_5 + \rho B_6 \big) (\rho)V,
     \end{equation}
where $\rho$ is the constant length of the parallel vector field $V$. As consequences, we have
\begin{itemize}
  \item [$\bullet$] $\tau _h (V)$ and $\tau _v (V)$ are homothetic to $V$ with constant homothety factors $c_h:= \big(2A_2 +mA_4  +\rho A_5 \big)(\rho)$ and $c_v:=\big(2B_3 +mB_5 + \rho B_6 \big) (\rho)$;
  \item [$\bullet$] $V:(M,g)\longrightarrow (TM,G)$ is a harmonic map if and only if $c_h =c_v =0$;
  \item [$\bullet$] $\tau _h (V)$ and $\tau _v (V)$ are parallel.
\end{itemize}

\begin{The}\label{gbitension_field}
Let $(M,g)$ be a Riemannian manifold, and $TM$ be equipped with a pseudo-Riemannian $g$-natural metric $G$. Let $V$ be a parallel vector field on $M$, with (constant) squared norm $\rho$. Then the bitension field of $V:(M,g) \rightarrow (TM,G)$ is given by $\tau_{2}(V)=(\tau_{2h}(V))^h \circ V +(\tau_{2v}(V))^v \circ V$, where
         \begin{eqnarray}
        \tau_{2h}(V) &=&\left\{(c_h^2 \rho -2 c_v) \rho[2A_2 +A_4 +\rho A_5](\rho) +2c_hc_v[C_2+ C_3+ C_5+ \rho C_6](\rho)\right. \label{hbitension}\\
        &&\left.+c_v^2[2E_1+ E_2+ \rho E_3](\rho) -2c_v \rho [2A_2 +A_4 +\rho A_5]^\prime(\rho)\right\}V,  \nonumber \\  \nonumber \\
          \tau_{2v}(V) &=&\left\{(c_h^2 \rho -2 c_v) \rho[2B_3 +B_5 +\rho B_6](\rho) +2c_hc_v[D_2+ D_3+ D_5+ \rho D_6](\rho)\right. \label{vbitension}\\
        &&\left.+c_v^2[2F_1+ F_2+ \rho F_3](\rho) -2c_v \rho [2B_3 +B_5 +\rho B_6]^\prime(\rho)\right\}V, \nonumber
         \end{eqnarray}
and $V:(M,g) \rightarrow (TM,G)$ defines a biharmonic map if and only if $\tau_{2h}(V)=\tau_{2v}(V)=0$.
         \end{The}

\begin{Proof}
Denote by $\widehat{\nabla}$ the induced connexion on $V^{-1}TTM$, and by $\widehat{\triangle}$ its associated rough laplacian. By (\ref{bitension_field}), we have:
         \begin{eqnarray}
                 \tau_{2}(V)&=&\widehat{\Delta}\tau(V)-\underset{i=1}{\overset{m}{\sum}}\bar{R}\big(\tau(V),V_{*}e_{i}\big)V_{*}e_{i} \nonumber \\
                 &=&-\underset{i=1}{\overset{m}{\sum}}\big(\widehat{\nabla}_{e_{i}}\widehat{\nabla}_{e_{i}}\tau(V) -\widehat{\nabla}_{\nabla_{e_{i}}e_{i}}\tau(V)\big)-\underset{i=1}{\overset{m}{\sum}}\bar{R}\big(\tau(V),V_{*}e_{i}\big)V_{*}e_{i}. \label{0}
                 \end{eqnarray}
We calculate each term of the sum separately using (\ref{cnx_exp}) and Proposition \ref{curvature_exp}. By definition of the induced connection $\nabla$, and taking into account the fact that $V_{*}X=X^{h}+(\nabla_{X}V)^{v}= X^h$, for every vector field $Y$, we have
\arraycolsep1.5pt
      \begin{eqnarray*}
     \widehat{\nabla}_{e_{i}}\widehat{\nabla}_{e_{i}}\tau(V) &=& \widehat{\nabla}_{e_{i}}\widehat{\nabla}_{e_{i}}\left((\tau_h(V))^h +(\tau_v(V))^v\right) \circ V \\
     &=& \left(\overline{\nabla}_{V_* e_i}\overline{\nabla}_{V_* e_i}\left((\tau_h(V))^h +(\tau_v(V))^v\right)\right) \circ V\\
     &=& \left(\overline{\nabla}_{e_{i}^{h}}\overline{\nabla}_{e_{i}^{h}}\left((\tau_h(V))^h +(\tau_v(V))^v\right)\right) \circ V.
      \end{eqnarray*}
Using Proposition \ref{lev-civ-con} and the fact that $\tau_h$ and $\tau_v$ are parallel, we deduce that
      \begin{eqnarray*}
     \widehat{\nabla}_{e_{i}}\widehat{\nabla}_{e_{i}}\tau(V) &=&  \left(\overline{\nabla}_{e_{i}^{h}} \Big\{h\big\{A\big(V;e_{i},\tau_{h}(V)\big) +C\big(V;e_{i},\tau_{v}(V)\big) \big\}\right.\\
     & &\left.+v\big\{B\big(V;e_{i},\tau_{h}(V)\big) +D\big(V;e_{i},\tau_{v}(V)\big)\big\}\Big\} \right)\circ V.
      \end{eqnarray*}
Using lemmas \ref{der-diff1} and \ref{der-F-T}, we have
      \begin{eqnarray*}
        \overline{\nabla}_{e_{i}^{h}} h\big\{A\big(V;e_{i},\tau_{h}(V)\big)\big\} &=& h\big\{\big(\nabla_{e_i} P^A_V\big)( e_{i},\tau_{h}(V))+ A\big(V;e_{i},A(V;e_i,\tau_{h}(V))\big)\big\} \\
        & & +v\big\{B\big(V;e_{i},A(V;e_i,\tau_{h}(V))\big)\big\}\\
         &=& h\big\{ A\big(V;e_{i},A(u;e_i,\tau_{h}(V))\big)\big\} +v\big\{B\big(V;e_{i},A(V;e_i,\tau_{h}(V))\big)\big\}.
      \end{eqnarray*}
In similar way, we calculate the other terms to obtain
      \begin{equation}\label{01}
        \begin{split}
     \widehat{\nabla}_{e_{i}}\widehat{\nabla}_{e_{i}}\tau(V) =&  h\Bigg\{ A\Big(V;e_{i},A\big(V;e_{i},\tau_{h}(V)\big)\Big)+C\Big(V;e_{i},B\big(V;e_{i},\tau_{h}(V)\big)\Big)\\
       &A\Big(V;e_{i},C\big(V;e_{i},\tau_{v}(V)\big)\Big)+C\Big(V;e_{i},D\big(V;e_{i},\tau_{v}(V)\big)\Big)\Bigg\} \\
       &+v\Bigg\{B\Big(V;e_{i},A\big(V;e_{i},\tau_{h}(V)\big)\Big) +D\Big(V;e_{i},B\big(V;e_{i},\tau_{h}(V)\big)\Big)\\
       &+B\Big(V;e_{i},C\big(V;e_{i},\tau_{v}(V)\big)\Big) +D\Big(V;e_{i},D\big(V;e_{i},\tau_{v}(V)\big)\Big)\Bigg\}.
         \end{split}
     \end{equation}

On the other hand, using Proposition \ref{lev-civ-con} and the fact that $\tau_h$ and $\tau_v$ are parallel, we have
     \begin{equation}\label{02}
        \begin{split}
          \widehat{\nabla}_{\nabla_{e_{i}}e_{i}}\tau(V) =& \widehat{\nabla}_{\nabla_{e_{i}}e_{i}}\left((\tau_h(V))^h +(\tau_v(V))^v\right) \circ V \\
          = & \Big(\overline{\nabla}_{(\nabla_{e_{i}}e_{i})^{h}} \left((\tau_h(V))^h +(\tau_v(V))^v\right)\Big) \circ V\\
         =&h\Bigg\{A\big(V;\nabla_{e_{i}}e_{i},\tau_{h}(V)\big)+ C\big(V;\nabla_{e_{i}}e_{i},\tau_{v}(V)\big) \Bigg\}\\ &+v\Bigg\{B\big(V;\nabla_{e_{i}}e_{i},\tau_{h}(V)\big)+D\big(V;\nabla_{e_{i}}e_{i},\tau_{v}(V)\big)\Bigg\}.
        \end{split}
     \end{equation}

Finally, using Proposition \ref{curvature_exp}, we have
      \begin{equation}\label{03}
        \begin{split}
     \bar{R}&\big(\tau(V),V_{*}e_{i}\big)V_{*}e_{i} = \bar{R}\big(\tau(V),e_{i}^{h}\big)e_{i}^{h}= \bar{R}\big(h\{\tau_{h}(V)\}+v\{\tau_{v}(V)\},e_{i}^{h}\big)e_{i}^{h}\\
     =&h\Bigg\{A\Big(V;\tau_{h}(V),A\big(V;e_{i},e_{i}\big)\Big)-A\Big(V;e_{i},A\big(V;\tau_{h}(V),e_{i}\big)\Big) +C\Big(V;\tau_{h}(V),B\big(V;e_{i},e_{i}\big)\Big) \\
      &-C\Big(V;e_{i},B\big(V;\tau_{h}(V),e_{i}\big)\Big)-A\Big(V;e_{i},C\big(V;e_{i}\,\tau_{v}(V)\big)\Big)-C\Big(V;e_{i},D\big(V;e_{i}\,\tau_{v}(V)\big)\Big)\\
      &+C\Big(V;A\big(V;e_{i},e_{i}\big),\tau_{v}(V)\Big)+E\Big(V;\tau_{v}(V),B\big(V;e_{i},e_{i}\big)\Big) -d\Big(A_{(e_{i},e_{i})}\Big)\big(\tau_{v}(V)\big)\Bigg\} \\
      &+v\Bigg\{B\Big(V;\tau_{h}(V),A\big(V;e_{i},e_{i}\big)\Big)-B\Big(V;e_{i},A\big(V;\tau_{h}(V),e_{i}\big)\Big) +D\Big(V;\tau_{h}(V),B\big(V;e_{i},e_{i}\big)\Big) \\
      &-D\Big(V;e_{i},B\big(V;\tau_{h}(V),e_{i}\big)\Big)-B\Big(V;e_{i},C\big(V;e_{i}\,\tau_{v}(V)\big)\Big)-D\Big(V;e_{i},D\big(V;e_{i}\,\tau_{v}(V)\big)\Big)\\
      &+D\Big(V;A\big(V;e_{i},e_{i}\big),\tau_{v}(V)\Big)+F\Big(V;\tau_{v}(V),B\big(V;e_{i},e_{i}\big)\Big)-d\Big(B_{(e_{i},e_{i})}\Big)\big(\tau_{v}(V)\big)\Bigg\}
         \end{split}
     \end{equation}

Substituting from \eqref{01}-\eqref{03} into \eqref{0}, we obtain
     \begin{equation*}
       \begin{split}
         \tau_{2}(V) = & \Big\{A\Big(V;\tau_{h}(V),\tau_{h}(V)\Big)+2C\Big(V;\tau_{h}(V),\tau_{v}(V)\Big)+E\Big(V;\tau_{v}(V),\tau_{v}(V)\Big)\\
           & -\underset{i=1}{\overset{m}{\sum}}d\Big(A_{(e_{i},e_{i})}\Big)\big(\tau_{v}(V)\big)\Big\}^h \circ V + \Big\{B\Big(V;\tau_{h}(V),\tau_{h}(V)\Big)+2D\Big(V;\tau_{h}(V),\tau_{v}(V)\Big) \\
           & +F\Big(V;\tau_{v}(V),\tau_{v}(V)\Big) -\underset{i=1}{\overset{m}{\sum}}d\Big(B_{(e_{i},e_{i})}\Big)\big(\tau_{v}(V)\big)\Big\}^v \circ V.
       \end{split}
     \end{equation*}
     Substituting from \eqref{tenspar1} into the last identity and using Proposition \ref{lev-civ-con} and Lemma \ref{der-diff2}, we obtain the result.
\end{Proof}


\subsection{Case of Kaluza-Klein type metrics}


           \begin{The}\label{biharm-KKtype}
                 Let $(M,g)$ be a Riemannian manifold and $V$ be a non zero parallel vector field on $M$ of (constant) squared norm $\rho$. Let the tangent bundle $TM$ of $M$ be equipped with a pseudo-Riemannian Kaluza-Klein type metric such that $(\beta_{1}+\beta_{3})(\rho)=(\beta_{1}+\beta_{3})'(\rho)=0$ (In particular, if $G$ is a Kaluza-Klein metric on $TM$). Then $V:(M,g) \rightarrow (TM,G)$ is a biharmonic map  if and only if one of the two statements is verified:
                       \begin{enumerate}
                       \item $\rho$ is a critical point of $\alpha_{1}+\alpha_{3}$.
                       \item $(\alpha_{1}+\alpha_{3})^{'}(\rho) \Big[2\rho \frac{\alpha_{1}^{'}}{\alpha_{1}}+\rho \frac{\beta_{1}-\alpha_{1}^{'}}{\phi_{1}}+\rho^{2} \frac{\alpha_{1} \beta_{1}^{'}-2\alpha_{1}^{'} \beta_{1}}{\phi_{1} \alpha_{1}}+2\rho\frac{\phi_{1}^{'}}{\phi_{1}} -1   \Big](\rho)-2\rho(\alpha_{1}+\alpha_{3})^{''}(\rho)=0 $.
                       \end{enumerate}
                \end{The}

\begin{Proof}
    If $\alpha_{2}=\beta_{2}=0$ and $(\beta_{1}+\beta_{3})(\rho)=(\beta_{1}+\beta_{3})'(\rho)=0$, then we have:
    \begin{equation*}
      B_{5}(\rho)=-\frac{(\alpha_{1}+\alpha_{3})^{'}}{\phi_{1}}(\rho), \quad  C_{3}(\rho)=-\frac{(\alpha_{1}+\alpha_{3})^{'}}{\alpha_{1}+\alpha_{3}}(\rho),
    \end{equation*}
    \begin{equation*}
      F_{1}(\rho)=\frac{\alpha_{1}^{'}}{\alpha_{1}}(\rho), \quad F_{2}(\rho)=\frac{\beta_{1}-\alpha_{1}^{'}}{\phi_{1}}(\rho), \quad F_{3}(\rho)=\frac{\alpha_{1} \beta_{1}^{'}-2\alpha_{1}^{'}\beta_{1} }{\phi_{1} \alpha_{1}}(\rho),
    \end{equation*}
    \begin{equation*}
      A_{2}(\rho)=A_{4}(\rho)=A_{5}(\rho)=B_{3}(\rho)=B_{6}(\rho)=C_{2}(\rho)=C_{5}(\rho)=C_{6}(\rho)=E_{1}(\rho)=E_{2}(\rho)=E_{3}(\rho)=0.
    \end{equation*}
    We deduce that $c_h=0$ and $c_v=-m\frac{(\alpha_{1}+\alpha_{3})^{'}}{\phi_{1}}(\rho)$.  Considering the results above, formulas (\ref{hbitension}) and (\ref{vbitension}) become:
    \begin{equation*}
      \begin{split}
        \tau_{2h}(V)= & 0, \\
         \tau_{2v}(V)= & m^{2}\frac{(\alpha_{1}+\alpha_{3})^{'}}{\phi_{1}^{2}}(\rho)\bigg((\alpha_{1}+\alpha_{3})^{'}(\rho) \Big[2\rho \frac{\alpha_{1}^{'}}{\alpha_{1}}+\rho \frac{\beta_{1}-\alpha_{1}^{'}}{\phi_{1}}+\rho^{2} \frac{\alpha_{1} \beta_{1}^{'}-2\alpha_{1}^{'} \beta_{1}}{\phi_{1} \alpha_{1}}+2\rho\frac{\phi_{1}^{'}}{\phi_{1}} -1   \Big](\rho)\\
         &\qquad\qquad\qquad-2\rho(\alpha_{1}+\alpha_{3})^{''}(\rho)\bigg)V.
      \end{split}
    \end{equation*}

    Since $\tau_{2h}(V)=0$ then $V:(M,g)\longrightarrow(TM,G)$ is a biharmonic map if and only if $\tau_{2v}(V)=0$, which gives the result.
\end{Proof}

\begin{Rq}\label{Rq-biharm}
  Let $G$ be the pseudo-Riemannian $g$-natural metric on $TM$ given by $\alpha_2=\beta_1=\beta_2=\beta_3=0$, $\alpha_1$ is constant and $(\alpha_1 +\alpha_3)(t)=e^{at}$, where $a>0$. Then $(\alpha_1+\alpha_3)^\prime$ does not vanish and
  $$(\alpha_{1}+\alpha_{3})^{'}(t) \Big[2t \frac{\alpha_{1}^{'}}{\alpha_{1}}+t \frac{\beta_{1}-\alpha_{1}^{'}}{\phi_{1}}+t^{2} \frac{\alpha_{1} \beta_{1}^{'}-2\alpha_{1}^{'} \beta_{1}}{\phi_{1} \alpha_{1}}+2t\frac{\phi_{1}^{'}}{\phi_{1}} -1   \Big](t)-2t(\alpha_{1}+\alpha_{3})^{''}(t)=-a(1+2at)e^{at}, $$
  which does not vanish on $\mathbb{R}^+$. We deduce that the two conditions of Theorem \ref{biharm-KKtype} can not be verified, and consequently every non-zero parallel vector field on $TM$ is not biharmonic.
\end{Rq}

    As corollary, we have the following result which characterizes parallel vector field on $M$ which are \textbf{proper biharmonic maps}, i.e. non-harmonic biharmonic maps.
    \begin{Cr}\label{Biharm-KKtype2}
     Let $(M,g)$ be a Riemannian manifold and $V$ be a parallel vector field on $M$ of (constant) squared norm $\rho$. Let the tangent bundle $TM$ of $M$ be equipped with a pseudo-Riemannian Kaluza-Klein type metric such that $(\beta_{1}+\beta_{3})(\rho)=(\beta_{1}+\beta_{3})'(\rho)=0$ (In particular, if $G$ is a Kaluza-Klein metric on $TM$).  Then the map $V:(M,g)\longrightarrow (TM,G)$ is proper biharmonic if and only if
    \begin{enumerate}
    \item $(\alpha_{1}+\alpha_{3})'(\rho)\neq 0 $ and
    \item  $(\alpha_{1}+\alpha_{3})^{'}(\rho) \Big[2\rho \frac{\alpha_{1}^{'}}{\alpha_{1}}+\rho \frac{\beta_{1}-\alpha_{1}^{'}}{\phi_{1}}+\rho^{2} \frac{\alpha_{1} \beta_{1}^{'}-2\alpha_{1}^{'} \beta_{1}}{\phi_{1} \alpha_{1}}+2\rho\frac{\phi_{1}^{'}}{\phi_{1}} -1   \Big](\rho)- 2\rho(\alpha_{1}+\alpha_{3})^{''}(\rho) =0$.
    \end{enumerate}
        \end{Cr}

        \begin{Exs}\label{Exs-prop-biharm}
        Let $(M,g)$ be a Riemannian manifold.
        \begin{enumerate}
        \item Let $G$ be the Kaluza-Klein metric on $TM$ such that
        \begin{itemize}
          \item $\beta_1=0$,
          \item $\alpha_1$ is a positive constant and
          \item $(\alpha_1+\alpha_3)(t)=\frac23 t^{3/2}+c$, for all $t \in \mathbb{R}_+$, $c>0$.
        \end{itemize}
        Then it is easy to see that $(\alpha_1+\alpha_3)^\prime$ doesn't vanish on $\mathbb{R}^*_+$ and the differential equation
        $$(\alpha_{1}+\alpha_{3})^{'}(t) \Big[2t \frac{\alpha_{1}^{'}}{\alpha_{1}}+t \frac{\beta_{1}-\alpha_{1}^{'}}{\phi_{1}}+t^{2} \frac{\alpha_{1} \beta_{1}^{'}-2\alpha_{1}^{'} \beta_{1}}{\phi_{1} \alpha_{1}}+2t\frac{\phi_{1}^{'}}{\phi_{1}} -1   \Big](t)- t(\alpha_{1}+\alpha_{3})^{''}(t)=0$$
        is satisfied on $\mathbb{R}^*_+$. We deduce that every non-zero parallel vector field on $M$ is a proper biharmonic map.
        \item Let $G$ be the Kaluza-Klein metric on $TM$ which satisfies
        \begin{itemize}
          \item $\alpha_{1}(t)=e^{\frac{t}{\rho}}$,
          \item $\beta_{1}(t)=\rho e^{\frac{t}{\rho}}-et$ and
          \item $(\alpha_1+\alpha_{3})(t)=k(e^{\frac{t}{\rho}}+c)$,
        \end{itemize}
        for all $t \in \mathbb{R}_+$, where $\rho$, $k$ and $c$ are constants such that $\rho \neq 0$, $k \neq 0$ and $c>-e$.
         Then it is easy to see that $G$ is non-degenerate. Furthermore, $G$ is Riemannian if $k>0$ and pseudo-Riemannian of signature $(n,n)$ if $k<0$. Obviously, $(\alpha_1+\alpha_3)^\prime$ doesn't vanish. It is also easy to check that the function $(\alpha_{1}+\alpha_{3})^{'}(t) \Big[2t \frac{\alpha_{1}^{'}}{\alpha_{1}}+t \frac{\beta_{1}-\alpha_{1}^{'}}{\phi_{1}}+t^{2} \frac{\alpha_{1} \beta_{1}^{'}-2\alpha_{1}^{'} \beta_{1}}{\phi_{1} \alpha_{1}}+2t\frac{\phi_{1}^{'}}{\phi_{1}} -1   \Big](t)- t(\alpha_{1}+\alpha_{3})^{''}(t)$, defined on $\mathbb{R}_+$, vanishes on $\rho$ and only on $\rho$. It follows that only parallel vector fields on $M$ of squared norm $\rho$ are (proper) biharmonic maps.
        \end{enumerate}
        \end{Exs}


\section{$G$-biharmonic vector fields}


    Let $(M,g)$ be a compact Riemannian manifold. In \cite{Ur.sect}, Markellos and Urakawa have proved that, equipping $TM$ with the Sasaki metric $g^{s}$, $V: (M,g) \longrightarrow (TM,g^{s})$ is a critical point of the bienergy functional restricted to $\mathfrak{X}(M)$ if and only if $V$ is a parallel vector field.\\
   Now, we shall determine the critical point condition for the bienergy functional restricted to $\mathfrak{X}(M)$, when $TM$ is equipped with an arbitrary pseudo-Riemannian $g$-natural metric. Such critical point will be called a \emph{$G$-biharmonic vector field}.

   Clearly, a $G$-biharmonic vector field is not necessarily a biharmonic map. The following result gives the characterization of $G$-biharmonic vector fields:

    \begin{The}\label{G-0}
    Let $(M,g)$ be a compact Riemannian manifold, and $TM$ be equipped with a pseudo-Riemannian $g$-natural metric $G$. A vector field $V$ is a $G$-biharmonic vector field if and only if
    \begin{equation}\label{G_cond}
    \alpha_{1}\tau_{2v}(V)+ \beta_{1} g\big(\tau_{2v}(V),V\big)V+\alpha_{2}\tau_{2h}(V)+\beta_{2} g\big(\tau_{2h}(V),V\big)V =0.
    \end{equation}
    In particular, a parallel vector field $V$ on $M$ is $G$-biharmonic if and only if
    \begin{equation}\label{G cond-par}
      \phi_{1}\tau_{2v}(V)=- \phi_{2}\tau_{2h}(V).
    \end{equation}
    \end{The}

\begin{Proof}
In \cite{Ji}, G. Jiang has derived the first variational formula of the bienergy functional $E_{2}$:
    $$\frac{d}{dt}\Big|_{t=0} E_{2}(V_{t})=-\int_{M}G\big(\tau_{2}(V),Y\big)dv_{g},$$
    for all $\{V_{t}\}\subset \mathfrak{X}(M)$ of $V$, and $Y$ is called the \textit{variational vector field} which satisfied:
    $$Y=\frac{\delta V_{t}}{\delta t}\Big|_{t=0},$$

    So $V$ is a $G$-biharmonic vector field i.e. a critical point of $E_{2\upharpoonleft \mathfrak{X}(M)}$ if and only if:
    $$\frac{d}{dt}\Big|_{t=0} E_{2}(V_{t})=0,$$
    for all $\{V_{t}\}\subset \mathfrak{X}(M)$ of $V$.

    But, as was remarked in \cite{Gilgilala}, any vertical vector field $Y^{v}$, section of the bundle $V^{-1}TTM$, there exists a variation $\{V_{t}\}_{t}\subset \mathfrak{X}(M)$ of $V$ such that
    $$Y^{v}=\frac{\delta V_{t}}{\delta t}\Big|_{t=0},$$

    Using the decomposition of $\tau_2(v)$ into horizontal and vertical components, we obtain

    \begin{equation*}
       \begin{split}
     \displaystyle \displaystyle \int_{M}G\big(\tau_{2}(V),Y^{v}\big)dv_{g}=&\displaystyle \displaystyle \int_{M}G\big(h\{ \tau_{2h}(V) \}+v\{ \tau_{2v}(V) \},Y^{v}\big)dv_{g}\\
     =&\displaystyle \displaystyle \int_{M} g\big(\alpha_{1}\tau_{2v}(V)+ \beta_{1} g\big(\tau_{2v}(V),V\big)V+\alpha_{2}\tau_{2h}(V)+\beta_{2} g\big(\tau_{2h}(V),V\big)V ,Y \big)dv_{g}.
       \end{split}
    \end{equation*}
    Hence $V$ is a $G$-biharmonic vector field if and only if \eqref{G_cond} is satisfied.

    Furthermore, if $V$ is parallel then $\tau_{2h}$ and $\tau_{2v}$ are proportional to $V$, and then \eqref{G_cond} is equivalent to \eqref{G cond-par}.
    \end{Proof}

    \begin{Rq}
    In the non compact case we can define the $G$-biharmonicity of $V$ by the condition \eqref{G_cond}, since it has a tensorial character.
    \end{Rq}

    Since the variations are through $\mathfrak{X}(M)\subset C^{\infty}(M,TM)$, if $V:(M,g)\longrightarrow (TM,G)$ is a biharmonic map then $V$ is a  $G$-biharmonic vector field $V$. The converse holds for parallel vector fields, for some pseudo-Riemannian $g$-natural metrics, as the following result shows

    \begin{Cr}
     Let $(M,g)$ be a Riemannian manifold, $V$ be a parallel vector field on $M$ and suppose that the tangent bundle $TM$ of $M$ is equipped by a pseudo-Riemannian $g$-natural metric $G$ which satisfies $\alpha_{2}(\rho)=\beta_{2}(\rho)=(\beta_{1}+\beta_{3})(\rho)=0$, where $\rho=||V||^{2}$.\par
     Suppose that we have one of the following conditions:
    \begin{enumerate}
    \item $m(\alpha_{1}+\alpha_{3})'(\rho)+\rho(\beta_{1}+\beta_{3})'(\rho)=0$ i.e. $V$ is a harmonic map;
    \item $[m(\alpha_{1}+\alpha_{3})'(\rho)+\rho(\beta_{1}+\beta_{3})'(\rho)]E_{1}(\rho)+\frac{\rho}{2}E_{3}(\rho)=mA'_{4}(\rho)+\rho A'_{5}(\rho)t$.
    \end{enumerate}
    Then  $V: (M,g) \longrightarrow (TM,G)$ is a biharmonic map if and only if $V$ is a $G$-biharmonic vector field.
    \end{Cr}
    \begin{Ex}\label{Ex-equiv}
      If $G$ is a Kaluza-Klein metric then $E_{1}=E_{3}=A'_{4}=A'_{5}=0$, and then a parallel vector field is $G$-biharmonic if and only if it is a biharmonic map..
    \end{Ex}

    On the other hand, we shall prove that if we consider some classes of pseudo-Riemannian $g$-natural metrics, then almost all parallel vector fields are $G$-biharmonic but not biharmonic. More precisely, suppose that $\alpha_1=\beta_1=\beta_2=\beta_3=0$, and $\alpha_1+\alpha_3$ does not vanish on $\mathbb{R}^+$. To ensure that $G$ is non-degenerate, $\alpha_2$ should not vanish on $\mathbb{R}^+$, and should have the same sign. In this case, we have
    \begin{itemize}
      \item $\alpha_1=\beta_1=\beta_2=\beta_3=\phi_1=0$, $\phi_1+\phi_3=\alpha_1 +\alpha_3$, $\phi=\alpha=-\alpha^2$.
      \item both $\alpha_1+\alpha_3$ and $\alpha_2$ do not vanish on $\mathbb{R}^+$.
      \item $G$ is of signature $(m,m)$. Indeed, for an orthonormal frame field $\{e_i;i=1,...,m\}$ on $(M,g)$, it suffices to take the orthonormal frame field $\{E_i;i=1,...,2m\}$ of $(TM,G)$ given by
          \begin{equation*}
            E_i=\frac{1}{\sqrt{|\alpha_1+\alpha_3|}}e_i^h, \qquad E_{m+i}=\frac{1}{\sqrt{|\alpha_1+\alpha_3|}}\Big[e_i^h -\frac{\alpha_1+\alpha_3}{\alpha_2}e_i^v\Big].
          \end{equation*}
    \end{itemize}
    To simplify notations, we consider the functions
    \begin{equation}\label{notations}
      \lambda:= \frac{(\alpha_1+\alpha_3)^\prime}{\alpha_2}, \qquad \mu:= \frac{\alpha_1+\alpha_3}{\alpha_2}, \qquad \nu:= \frac{\alpha_2^\prime}{\alpha_2}.
    \end{equation}
    Since we are looking for non biharmonic vector fields on $M$, we will suppose that $(\alpha_1+\alpha_3)^\prime$ does not vanish on $\mathbb{R}^*_+$. We deduce that $\lambda$ and $\mu$ don't vanish on $\mathbb{R}^*_+$. Furthermore, taking eventually $-\alpha_2$ instead of $\alpha_2$, we can suppose that $\lambda$ is positive on $\mathbb{R}^*_+$.

    Substituting from \eqref{notations} in Proposition \ref{lev-civ-con} and in the formulas of $c_h$ and $c_v$, we find
    \begin{itemize}
      \item $A_1=A_2=A_3=A_5=0$, $A_4= -\lambda$,
      \item $B_1=-1$, $B_3=B_4=B_6=0$, $B_5 =\lambda\mu$,
      \item $C_1=C_2=C_4=C_6=0$, $C_3=\nu$, $C_5=-\nu$,
      \item $D1=D_2=D_4=D_6=0$, $D_3=\mu^\prime$, $D_5=-\mu\nu$,
      \item $E_1=E_2=E_3=0$,
      \item $F_1=\nu$, $F_2=F_3=0$,
      \item $c_h=-m\lambda$, $c_v= m\lambda \mu$.
    \end{itemize}
    Substituting into \eqref{hbitension} and \eqref{vbitension}, we deduce that the horizontal and vertical components of the bitension of a parallel vector field $V$ om $M$ with $\rho:=\|V|^2$ are given by
    \begin{eqnarray}
      \tau_{2h}(V) &=& m \lambda(\rho)[2\mu(\lambda +\rho\lambda^\prime) -m\rho \lambda^2](\rho)V, \label{G-1} \\
      \tau_{2v}(V) &=& m \lambda(\rho)\mu(\rho)\{m\rho^2\lambda^2 -2\rho[\lambda\mu +(\lambda\mu)^\prime] +m\lambda(\rho)[3\mu\nu -2\mu^\prime]\}(\rho)V. \label{G-2}
    \end{eqnarray}
    From Theorem \ref{G-0}, $V$ is $G$-biharmonic if and only if $\tau_{2h}(V)=0$. Hence $V$ is non biharmonic $G$-biharmonic vector field if and only if $\tau_{2h}(V)=0$ and $\tau_{2v}(V)\neq0$, i.e. if and only if the two following conditions hold
    \begin{eqnarray}
      & & [2\mu(\lambda +\rho\lambda^\prime) -m\rho \lambda^2](\rho)=0, \label{G-3} \\
      & & \{m\rho^2\lambda^2 -2\rho[\lambda\mu +(\lambda\mu)^\prime] +m\lambda(\rho)[3\mu\nu -2\mu^\prime]\}(\rho) \neq 0. \label{G-4}
    \end{eqnarray}
    Now, we shall consider the differential equation $2\mu(\lambda +t\lambda^\prime) -mt \lambda^2=0$, which is equivalent to
    \begin{equation}\label{G-5}
      \frac{\lambda}{\mu}= \frac{2}{m}\Big(\frac{1}{t} +\frac{\lambda^\prime}{\lambda}\Big)
    \end{equation}
     on $\mathbb{R}^*_+$. Remarking that $\frac{\lambda}{\mu}= \frac{(\alpha_1 +\alpha_3)^\prime}{\alpha_1 +\alpha_3}$ and integrating, we find that
     \begin{equation}\label{G-6}
       \alpha_1 +\alpha_3=K (t\lambda)^{\frac2m},
     \end{equation}
     for some non zero constant $K$. Substituting from \eqref{G-5} and \eqref{G-6} into $\mu=\frac{\alpha_1+\alpha_3}{\alpha_2}$, we obtain
     \begin{equation}\label{G-7}
       \alpha_2=\frac{2K}{m}(\lambda +t\lambda^\prime) t^{\frac2m-1}\lambda^{\frac2m-2},
     \end{equation}
     on $\mathbb{R}^*_+$. To avoid the non regularity problem of $\alpha_2$ at 0, we need to have $\lim\limits_{t \rightarrow 0} \frac{\lambda(t)}{t}$ is a constant. In this case, we will have $\alpha_2(0)=0$, which is in contradiction with the fact that $G$ is non-degenerate at the zero section. To solve this problem, we consider a functions $f_\eta \in C^\infty(\mathbb{R}^+)$ such that
     \begin{itemize}
       \item $0\leq f_\eta(t) \leq 1$, for all $t \in \mathbb{R}^+$,
       \item $f_\eta(t)=0$, for $t \geq \eta$,
       \item $f_\eta(t)=1$, for $t \leq \frac{\eta}{2}$,
     \end{itemize}
     and we put
     \begin{equation}\label{G-8}
       \left\{
       \begin{split}
          &\alpha_1+\alpha_3 = K (t\lambda)^{\frac2m} +\textup{sgn}(K)f_\eta, \\
          & \alpha_2=\frac{2K}{m}(\lambda +t\lambda^\prime) t^{\frac2m-1}\lambda^{\frac2m-2} +\textup{sgn}(K(\lambda +t\lambda^\prime))f_\eta.
       \end{split}
       \right.
     \end{equation}
     Hence we have
     \begin{itemize}
       \item $\alpha_1+\alpha_3$ and $\alpha_2$ are exactly those of \eqref{G-6} and \eqref{G-7} on $[\eta,+\infty[$, respectively, and satisfy the condition \eqref{G-3} on $[\eta,+\infty[$. We deduce that, for $\rho \in [\eta,+\infty[$, a parallel vector of norm $\rho$ is $G$-biharmonic.
       \item $\alpha_2(0) \neq 0$ and then $G$ is non-degenerate everywhere.
       \item $\alpha_1+\alpha_3$ and $\alpha_2$ don't vanish on $\mathbb{R}^+$.
     \end{itemize}
     To complete our study, we shall check if condition \eqref{G-4} is satisfied for $\rho \in [\eta,+\infty[$. Using \eqref{G-8} on $[\eta,+\infty[$, we can express $\mu$ and $\nu$ in terms of  $\lambda$ and its derivatives as
     \begin{eqnarray*}
       \mu &=& \frac{mt\lambda^2}{2(\lambda +t\lambda^\prime)}, \\
       \nu &=&  \frac{1}{t\lambda(\lambda +t\lambda^\prime)}\Big[t^2\lambda\lambda^{\prime\prime} +\Big(\frac4m -1\Big)t\lambda\lambda^\prime +\Big(\frac2m -2\Big)t^2 \big(\lambda^\prime\big)^2 +\Big(\frac2m -1\Big)\lambda^2\Big].
     \end{eqnarray*}
     We deduce then that condition \eqref{G-4} at $\rho \in [\eta,+\infty[$ is equivalent to the following condition:
     \begin{equation}\label{G-9}
       \begin{split}
         0\neq  & \left\{t^2(\lambda +t\lambda^\prime)^2 -t\Big[((1+t)\lambda^3 +3\lambda^2 \lambda^\prime) -t\lambda^3(\lambda +t\lambda^\prime)^\prime\Big] +3m\lambda^2\Big[t^2\lambda^{\prime\prime}+\Big(\frac4m -1\Big)t\lambda\lambda^\prime \right. \\
           & \left.+\Big(\frac2m -2\Big)t^2 \big(\lambda^\prime\big)^2 +\Big(\frac2m -1\Big)\lambda^2\Big] -m\lambda (\lambda^2 +2t\lambda\lambda^\prime)(\lambda +t\lambda^\prime) +mt\lambda^3(\lambda +t\lambda^\prime)^\prime\right\}(\rho).
       \end{split}
     \end{equation}
     To summarize the previous discussion, we have
     \begin{Prop}\label{G-biharm-not-biharm}
       Let $(M,g)$ be a Riemannian manifold. Fix $\eta>0$ and suppose that the tangent bundle $TM$ of $M$ is equipped by a pseudo-Riemannian $g$-natural metric $G$ which satisfies $\alpha_1=\beta_1=\beta_2=\beta_3=0$ and \eqref{G-8}, where $\lambda$ is a positive smooth function on $\mathbb{R}^*_+$ such that $\lim\limits_{t \rightarrow 0} \frac{\lambda(t)}{t}$ is a constant. Then a parallel vector field on $M$ such that $\rho:=||V||^{2} \in [\eta,+\infty[$ is non-biharmonic $G$-biharmonic if and only if \eqref{G-9} is satisfied at $\rho$.
     \end{Prop}

     \begin{Ex}\label{Ex-G-biharm-not-biharm}
       We consider the function $\lambda$ defined by $\lambda(t) =ate^{bt}$, for all $t \in \mathbb{R}^+$, where $a>0$ and $b$ are real constants. A long but routine calculation shows that \eqref{G-9} is equivalent to the fact that $\rho$ is not a solution of the following polynomial equation
       \begin{equation*}
         [(2m-11)b-1]bt^3 +[2mb^2 +(5m-21)b +2]t^2 +[7mb +17m -20]t -6m = 0,
       \end{equation*}
       which has at most three solutions. We deduce that, apart from at most three values of $\rho \in [\eta,+\infty[$, all parallel vector fields of norm $\sqrt{\rho}$ are non-biharmonic $G$-biharmonic.
     \end{Ex}
     \begin{Rqs}
       \begin{enumerate}
         \item Taking eventually different values of $b$ in the previous example, for any value $\rho \in [\eta,+\infty[$, there is a pseudo-Riemannian $g$-natural metric on $TM$ for which any parallel vector field on $M$ such that $||V||^{2}=\rho$ is non-biharmonic $G$-biharmonic.
         \item $\eta$ being arbitrary, for any value $\rho >0$, there is a pseudo-Riemannian $g$-natural metric on $TM$ for which any parallel vector field on $M$ such that $||V||^{2}=\rho$ is non-biharmonic $G$-biharmonic.
         \item It would be interesting to study the two following (not easy) problems:
         \begin{itemize}
           \item Are there examples of non biharmonic $G$-biharmonic vector fields when the tangent bundle is equipped with an appropriate \textbf{Riemannian} $g$-natural metric?
           \item Can we find pseudo-Riemannian $g$-natural metrics for which there are non parallel biharmonic (resp. $G$-biharmonic) vector fields?
         \end{itemize}

       \end{enumerate}
     \end{Rqs}


\section{Biharmonic unit vector fields}


\subsection{The main theorem}


    \textbf{Notations:}
    Let $(M,g)$ be a Riemannian manifold and $U\in \mathfrak{X}_{1}(M)$. Let $\{e_{i}\}^{m}_{i=1}$ be an orthonormal frame field on $M$. Let us denote by
    \begin{itemize}
      \item [$\bullet$] $S(U)=-\sum_{i=1}^{m}R(\nabla_{e_{i}}U,U)e_{i}$,
      \item [$\bullet$] $Q$ the Ricci operator associated to $g$,
      \item [$\bullet$] $\textup{div}$ the divergence operator associated to $g$,
      \item [$\bullet$] $\mathcal{A}$ the operator on unit vector fields defined by: $\mathcal{A}U=QU-g(QU,U)U$.
    \end{itemize}

   Let $(M,g)$ be a Riemannian manifold and $T_{1}M$ its unit tangent bundle. We call \emph{$g$-natural metrics} on $T_{1}M$ the restrictions of $g$-natural metrics of $TM$ to its hypersurface $T_{1}M$. These metrics possess a simpler form. Indeed, for any $g$-natural metric  $\tilde{G}$ on $T_{1}M$, there are four real constants $a$, $b$, $c$ and $d$, such that
         \arraycolsep1.5pt
            \begin{equation}\label{exp-$g$-nat}
             \left\lbrace
            \begin{array}{rcl}
             \tilde{G}_{(x,u)}(X^{h},Y^{h})&=&(a+c)g_{x}(X,Y)+dg_{x}(X,u)g_{x}(Y,u),\\
             \tilde{G}_{(x,u)}(X^{h},Y^{t})&=&b g_{x}(X,Y),\\
             \tilde{G}_{(x,u)}(X^{t},Y^{t})&=&ag_{x}(X,Y)-\frac{\phi}{a+c+d}g_{x}(X,u)g_{x}(Y,u),
            \end{array}
            \right.
            \end{equation}
for any $(x,u) \in T_{1}M$ and $X$, $Y\in M_{x}$ (cf. \cite{Abb-Kow1}), where $X^t$ is the tangential lift to $(x,u)$ of $X$ given by
\begin{equation}\label{tang1}
    X^{t} = \frac{b}{a+c+d} g_x(X,u) \,u^h + [X -  g_x(X,u)u]^v.
\end{equation}

Furthermore, $\tilde{G}$ is non-degenerate (resp. Riemannian) if and only if $\alpha:=a(a+c)-b^{2} \neq 0$ and $\phi=a(a+c+d)-b^{2} \neq 0$ (resp. $a>0$, $\phi >0$ and $\alpha>0$).

In particular, the {\em Sasaki metric} on $T_1 M$ corresponds to the case where $a=1$ and $b=c=d=0$; {\em Kaluza-Klein metrics} are obtained when
$b=d=0$; {\em metrics of Kaluza-Klein type} are given by the case $b=0$.

From now on, let $(M,g)$ be a compact Riemannian manifold and suppose its unit tangent sphere bundle $T_{1}M$ is equipped by a pseudo-Riemannian $g$-natural metric $\tilde{G}$.\par
The \emph{bienergy} functional of a unit vector field $U$, regarded as a map of $(M,g)$ into $(T_{1}M,\tilde{G})$, is the mapping defined by:
      \begin{equation*}
        \begin{split}
          E_{2}: C^{\infty}(M,T_{1}M) & \longrightarrow [0,+\infty[ \\
           U & \longmapsto E_{2}(U):=\frac{1}{2}\int_{M}\|\tilde{\tau}(U)\|^{2}v_{g},
        \end{split}
      \end{equation*}
      $\tilde{\tau}(U)$ being the tension field of the map $U:(M,g) \longrightarrow (T_{1}M,\tilde{G})$, which is given by: $\tilde{\tau}(U)=(\tilde{\tau}_{h}(U))^h+ (\tilde{\tau}_{v}(U))^v$, where (cf. \cite{Ab.unit})
       \begin{eqnarray}
                \tilde{\tau}_{h}(U)&=&\frac{ab}{\alpha}QU-\frac{a^{2}}{\alpha}S(U)-\frac{ad}{\alpha}\nabla_{U}U-\frac{b(ad+b^{2})}{\alpha \varphi}g(QU,U)U-\frac{b}{\varphi}g(\overline{\Delta}U,U)U  \label{11}\\
                &&+\frac{d}{\varphi}\textup{div}(U)U+\frac{a(ad+b^{2})}{\alpha \varphi}g(S(U),U)U \nonumber \\ \nonumber \\
              \tilde{\tau}_{v}(U)&=&-\overline{\Delta}U-\frac{b^{2}}{\alpha}QU+\frac{ab}{\alpha}S(U) +\frac{bd}{\alpha}\nabla_{U}U+\frac{b^{2}}{\alpha}g(QU,U)U+g(\overline{\Delta}U,U)U  \label{12}\\
              &&-\frac{ab}{\alpha}g(S(U),U)U, \nonumber
       \end{eqnarray}
         with  $\varphi=a+c+d.$

A unit vector field $U$ is said to be a \emph{$\tilde{G}$-biharmonic unit vector field} if $U$ is a critical point of $E_{2 \upharpoonright \mathfrak{X}_{1}(M)}$, i.e. for all $C^{\infty}$ 1-parameter variations $\{U_{t}\}_{t}$ of $U$ in $\mathfrak{X}^{1}(M)$, $(|t|<\epsilon)$, we have:
           $$\frac{d}{dt}\Big|_{t=0} E_{2}(U_{t})=0.\label{egale0}$$

    \begin{The}\label{buvcond}
    Let $(M,g)$ be a compact Riemannian manifold and $T_{1}M$ its unit tangent bundle equipped with an arbitrary pseudo-Riemannian $g$-natural metric $\tilde{G}$. Let $U \in \mathfrak{X}^{1}(M)$. Then $U$ is a $\tilde{G}$-biharmonic unit vector field if and only if the following holds
    $$T(U)=0,$$ where
     \begin{equation}\label{buvcond0}
     \begin{array}{rcl}
     &&T(U)=\underset{i=1}{\overset{m}{\sum}}\Big\{-a\big[2R(e_{i},\tilde{\tau}_{h}(U))\nabla_{e_{i}}U +(\nabla_{e_{i}}R)(e_{i},\tilde{\tau}_{h}(U))U+R(e_{i},\nabla_{e_{i}}\tilde{\tau}_{h}(U))U\big]\\\\

       &&-dg\big(\nabla_{e_{i}}U,\tilde{\tau}_{h}(U)\big)e_{i}\Big\}+d\:\nabla_{U}\tilde{\tau}_{h}(U)+b\:Q\big(\tilde{\tau}_{h}(U)\big)\\\\
        &&+\big[-\frac{bd}{\varphi}g(QU,U)+\frac{bd}{\varphi}g(\overline{\Delta}U,U)+(a+c)\frac{d}{\varphi}\textup{div}(U)+\frac{ad}{\varphi}g(S(U),U) +dg(\tilde{\tau}_{h}(U),U)\big]\tilde{\tau}_{h}(U)\\\\
        &&+\frac{a^{2}d}{\alpha} g(\tilde{\tau}_{h}(U),U)S(U)-\frac{abd}{\alpha} g(\tilde{\tau}_{h}(U),U)QU-a\overline{\Delta}\tilde{\tau}_{v}(U)-b\overline{\Delta}\tilde{\tau}_{h}(U)\\\\

        &&+\big[\frac{b^{2}}{\varphi}g(QU,U)+\frac{\alpha +ad}{\varphi}g(\overline{\Delta}U,U)-\frac{ab}{\varphi}g(S(U),U)+\frac{bd}{\varphi}\textup{div}(U)\big] \tilde{\tau}_{v}(U)   \\\\
           &&-d\: \textup{grad}[g(\tilde{\tau}_{h}(U),U)]+\frac{ad^{2}}{\alpha}g\big(\tilde{\tau}_{h}(U),U\big)\nabla_{U}U.
     \end{array}
     \end{equation}
    \end{The}

    \begin{Proof}
    Let $U$ be a unit vector field on $(M,g)$ and $I=]-\epsilon,\epsilon[$, $(\epsilon >0)$. For $t \in I$, we denote by $i_{t}:M \longrightarrow M \times I$, $p\longrightarrow(p,t)$, the canonical inclusion. We consider $C^{\infty}$-variations $V:M \times I \longrightarrow T_{1}M$ of $U$ within $\mathfrak{X}_{1}(M)$, i.e. for all $t \in I$, the mappings $V_{t}=V\circ i_{t}$ are in fact unit vector fields and $V_{0}=U$. We choose $\{e_{i}\}_{i=1}^{m}$ a local othonormal frame field of $(M,g)$ and we extend each $e_{i}$ (resp. $\frac{d}{dx}\in \mathfrak{X}(I)$) to $M \times I$, denoted by $E_{i}$ (resp. $\frac{d}{dt}$). If we denote by $D$ the Levi-Civita connection  and $R^{D}$ the Riemannian curvature of the Riemannian product manifold $M \times I$, then using the second Bianchi identity for the last relation, we get
    \begin{equation}\label{10}
     [E_{i},\frac{d}{dt}]=0, \qquad D_{\frac{d}{dt}}E_{i}=0, \qquad R^{D}(TM,TI)=0,  \qquad (D_{\frac{d}{dt}}R^{D})(D_{E_{i}}V,V)E_{i}=0,
    \end{equation}
    for all $1 \leq i \leq m$. We set $\tilde{\tau}_h:= \tilde{\tau}_{h}(V)$ and $\tilde{\tau}_v:= \tilde{\tau}_{v}(V)$. It is easy to see that  $\tilde{\tau}_{h}(V_t)=\tilde{\tau}_{h} \circ i_{t}$ and $\tilde{\tau}_{v}(V_t)=\tilde{\tau}_{v} \circ i_{t}$.

    We have to calculate the first variational formulae of $E_{2 \upharpoonright \mathfrak{X}^{1}(M)}$. By definition of $\tilde{G}$, we have
    \begin{eqnarray}
    E_{2}(V_{t})
       &=&\frac{a+c}{2}\displaystyle \int_{M}g\big(\tilde{\tau}_{h},\tilde{\tau}_{h}\big)\circ i_{t}dv_{g}+\frac{d}{2}\displaystyle \int_{M}g\big(\tilde{\tau}_{h},V\big)^{2}\circ i_{t}dv_{g}+b\displaystyle \int_{M}g\big(\tilde{\tau}_{h},\tilde{\tau}_{v}\big)\circ i_{t}dv_{g}\\
          &&+\frac{a}{2}\displaystyle \int_{M}g\big(\tilde{\tau}_{v},\tilde{\tau}_{v}\big)\circ i_{t}dv_{g}.\nonumber
    \end{eqnarray}

    Then we get
    \begin{eqnarray}
    D_{\frac{d}{dt}}E_{2}(V_{t})
    &=&(a+c)\displaystyle \int_{M}g(D_{\frac{d}{dt}}\tilde{\tau}_{h},\tilde{\tau}_{h})\circ i_{t}dv_{g}+d\displaystyle \int_{M}g(\tilde{\tau}_{h},V)g(D_{\frac{d}{dt}}\tilde{\tau}_{h},V)\circ i_{t}dv_{g}\\ \nonumber
    &&+d\displaystyle \int_{M}g(\tilde{\tau}_{h},V)g(\tilde{\tau}_{h},D_{\frac{d}{dt}}V)\circ i_{t}dv_{g}+b\displaystyle \int_{M}g(D_{\frac{d}{dt}}\tilde{\tau}_{h},\tilde{\tau}_{v})\circ i_{t}dv_{g}\\ \nonumber
    &&+b\displaystyle \int_{M}g(D_{\frac{d}{dt}}\tilde{\tau}_{v},\tilde{\tau}_{h})\circ i_{t}dv_{g}+a\displaystyle \int_{M}g(D_{\frac{d}{dt}}\tilde{\tau}_{v},\tilde{\tau}_{v})\circ i_{t}dv_{g}.
    \end{eqnarray}
    By long computation, using \eqref{11}, \eqref{12} and \eqref{10}, we get $T(U)=0$.
    \end{Proof}

    \begin{Rq}
      In the non-compact case, we can define the biharmonicity of unit vector fields by the condition $T(U)=0$, where $T$ is given by (\ref{buvcond0}), since it has a tensorial character.
        \end{Rq}

    \begin{Prop}
    Let $(M,g)$ be a Riemannian manifold and $T_{1}M$ be its unit tangent bundle equipped with an arbitrary pseudo-Riemannian $g$-natural metric $\tilde{G}$. Let $U$ be a unit vector field on $M$. If $U:(M,g) \longrightarrow (T_{1}M,\tilde{G})$ is a harmonic map then $U$ is a $\tilde{G}$-biharmonic unit vector field.
    \end{Prop}

    Note that the converse does not hold. $\tilde{G}$-biharmonic unit vector field which are not harmonic maps $(M,g) \rightarrow (T_1M,\tilde{G})$ are called a \emph{proper $\tilde{G}$-biharmonic unit vector field}.


    \subsection{Examples of $\tilde{G}$-biharmonic unit vector fields}


\begin{Ex}\label{Ex-hyperbolic}

Consider the hyperbolic space $\mathbb{H}^{n}$ $(n>1)$ of constant negative sectional curvature $-k^{2}$, that is, $\mathbb{H}^{n}=(\mathbb{R}^{n}_{+},g)$, where $\mathbb{R}^{n}_{+}=\{(y_{1},...,y_{n}) \in \mathbb{R}^{n}: y_{n}>0\}$ and $$g=\frac{1}{k^{2}y_{n}^{2}}(dy_{1}\otimes dy_{1}+...+dy_{n}\otimes dy_{n})$$
Vector fields $e_{i}=ky_{n}\frac{\partial}{\partial y_{i}}$ for $i=1,...,n$ provide an orthonormal frame on $\mathbb{H}^{n}$. Put $V=e_{n}$. A standard calculation shows that covariant derivatives of $e_{i}$ can be completly descrined as follows:
\begin{equation}\label{hypercarac}
\nabla_{e_{i}}e_{j}=k\delta_{ij}V, \phantom{XXXX}  \nabla_{e_{i}}V=-ke_{i}, \phantom{XXXX} \nabla_{V}e_{i}=0, \phantom{XXXX} \nabla_{V}V=0,
\end{equation}

for all $i, j<n$. In particular, (\ref{hypercarac}) easily implies:
\begin{equation}\label{hyper_formulas}
 ||\nabla V||^{2}=(n-1)k^{2},\phantom{X}  \overline{\Delta}V=-tr \nabla^{2}{V}=||\nabla V||^{2} V=(n-1)k^{2}V=QV,\phantom{X}
div(V)=(1-n)k,
\end{equation}
and $$ \phantom{X} S(V)=k^{2}div(V)V=k^{3}(1-n)V.$$
The vertical and horizontal parts of the tension field become:

\begin{equation}\label{hyper_tension_parts}
 \tilde{\tau}_{v}(V)=0, \phantom{XXXX}
 \tilde{\tau}_{h}(V)=\frac{k}{\varphi}(1-n)(d-ak^{2})V.
\end{equation}
Using (\ref{hypercarac}), (\ref{hyper_formulas}) and (\ref{hyper_tension_parts}), the condition $T(V)=0$ of the biharmonicity of a unit vector field $V$ on $M$ becomes:
 \begin{equation}\label{hyper-biharm-cond}
            \begin{array}{rcl}
            &&\frac{k^{2}}{\varphi}(1-n)^{2}(d-ak^{2})\Big\{(-2a+\frac{a^{2}d}{\alpha}  )k^{2}+\frac{abd}{\alpha} k+d\Big\}V=0.
            \end{array}
            \end{equation}
 It follows that
 \begin{itemize}
   \item either $d=ak^2$ and, in this case, $V$ is a harmonic map $(M,g) \rightarrow (T_1M,\tilde{G})$;
   \item or $d \neq ak^2$. In this case, $d \neq 0$ since if not we have $a \neq 0$ and \eqref{hyper-biharm-cond} gives then $ak^2=0$, which is a contradiction. On the other hand, $a \neq 0$ since if not \eqref{hyper-biharm-cond} yields $\alpha d=0$, which is a contradiction. Taking into account these two last facts, \eqref{hyper-biharm-cond} is equivalent to
       \begin{equation}\label{hyper-biharm-b}
         b=\left(\frac{2\alpha}{d} -a\right)k -\frac{\alpha}{ak}.
       \end{equation}
 \end{itemize}
 We conclude then the following
 \begin{Prop}
   Let $\mathbb{H}^{n}$ $(n>1)$ be the hyperbolic space of constant negative sectional curvature $-k^{2}$ and $\tilde{G}$ be an arbitrary pseudo-Riemannian $g$-natural metric on its unit tangent bundle. Then the vector field $V=e_{n}$ is a $\tilde{G}$-biharmonic unit vector field if and only if one of the two following statements is verified:
\begin{enumerate}
\item $d=ak^{2}$. In this case, $V$ is a harmonic map $(M,g) \rightarrow (T_1M,\tilde{G})$.
\item $a \neq 0$, $d \neq 0$, $d \neq ak^2$ and \eqref{hyper-biharm-b}. In this case, $V$ is a proper $\tilde{G}$-biharmonic unit vector field.
\end{enumerate}
 \end{Prop}

\begin{Cr}
Let $\mathbb{H}^{n}$ $(n>1)$ be the hyperbolic space of constant negative sectional curvature $-k^{2}$. Then there is a three-parameter family of pseudo-Riemannian $g$-natural metrics on the unit tangent bundle with respect to which the vector field $V=e_{n}$ is a proper $\tilde{G}$-biharmonic unit vector field.
\end{Cr}
\end{Ex}

\begin{Ex}\label{Ex-Solvable}

Consider the solvable Lie group \textit{Sol}$_{3}$ as ${\mathbb{R}}^{3}$ equipped with the metric  $g_{sol}=e^{2z}dx^{2}+e^{-2z}dy^{2}+dz^2$ and the orthonormal frame field on it \cite{sol3}:
\begin{equation}
e_{1}=e^{-z}\frac{\partial}{\partial x}, \phantom {XXX}e_{2}=e^{z}\frac{\partial}{\partial y},\phantom {XXX}e_{3}=\frac{\partial}{\partial z}.
\end{equation}
With respect to this orthonormal frame, the Lie brackets and the Levi-Civita connection can be easily computed as:
\begin{center}
\begin{tabular} {ccc}
$[e_{1},e_{2}]=0$,&$[e_{2},e_{3}]=-e_{2}$,&$[e_{1},e_{3}]=e_{1}$,\\\\
$\nabla_{e_{1}}e_{1}=-e_{3} $,&$\nabla_{e_{1}}e_{2}=0$,&$\nabla_{e_{1}}e_{3}=e_{1}$,\\\\
$\nabla_{e_{2}}e_{1}=0 $,&$\nabla_{e_{2}}e_{2}=e_{3}$,&$\nabla_{e_{2}}e_{3}=-e_{2}$,\\\\
$\nabla_{e_{3}}e_{1}=0 $,&$\nabla_{e_{3}}e_{2}=0$,&$\nabla_{e_{3}}e_{3}=0$.\\\\
\end{tabular}
\end{center}
A further computation gives:
\begin{center}
\begin{tabular}{ccc}
$R(e_{1},e_{2})e_{1}=-e_{2}$,&$R(e_{1},e_{3})e_{1}=e_{3}$,&$R(e_{1},e_{2})e_{2}=e_{1}$,\\\\
$R(e_{2},e_{3})e_{2}=e_{3}$,&$R(e_{1},e_{3})e_{3}=-e_{1}$,&$R(e_{2},e_{3})e_{3}=-e_{2}$.\\\\
\end{tabular}
\end{center}
and the terms in $\ref{buvcond}$ become:
\begin{center}
\begin{tabular}{c}
$\underset{i=1}{\overset{m}{\sum}}\Big\{2R(e_{i},e_{1})\nabla_{e_{i}}e_{1}+(\nabla_{e_{i}}R)(e_{i},e_{1})e_{1}+R(e_{i},\nabla_{e_{i}}e_{1})e_{1}\Big\}=e_{3}$,\\\\
$\underset{i=1}{\overset{m}{\sum}}\Big\{2R(e_{i},e_{2})\nabla_{e_{i}}e_{2}+(\nabla_{e_{i}}R)(e_{i},e_{2})e_{2}+R(e_{i},\nabla_{e_{i}}e_{2})e_{2}\Big\}=0$,\\\\
$\underset{i=1}{\overset{m}{\sum}}\Big\{2R(e_{i},e_{3})\nabla_{e_{i}}e_{3}+(\nabla_{e_{i}}R)(e_{i},e_{3})e_{2}+R(e_{i},\nabla_{e_{i}}e_{3})e_{3}\Big\}=0$.
\end{tabular}
\end{center}
\begin{center}
\begin{tabular}{ccc}
$\overline{\Delta}e_{1}=e_{1}$,&$\overline{\Delta}e_{2}=e_{2}$,&$\overline{\Delta}e_{3}=2e_{3}$,\\\\
$Q(e_{1})=0$,&$Q(e_{2})=0$,&$Q(e_{3})=-2e_{3}$,\\\\
$div(e_{1})=0$,&$div(e_{2})=0$,&$div(e_{3})=0$,\\\\
$S(e_{1})=-e_{3}$,&$S(e_{2})=e_{3}$,&$S(e_{3})=0$,\\\\
$\tilde{\tau}_{h}(e_{1})=(\frac{a(a+d)}{\alpha}-\frac{b}{\varphi})e_{3}$,&$\tilde{\tau}_{h}(e_{2})=-\frac{a}{\alpha}(a+d)e_{3}-\frac{b}{\varphi}e_{2}$,&$\tilde{\tau}_{h}(e_{3}) =-\frac{4b}{\varphi}e_{3}$,\\\\
$\tilde{\tau}_{v}(e_{1})=-\frac{b}{\alpha}(a+d)e_{3}$,&$\tilde{\tau}_{v}(e_{2})=\frac{b}{\alpha}(a+d)e_{3}$,&$\tilde{\tau}_{v}(e_{3})=0$.
\end{tabular}
\end{center}
Using the formulas above, the condition of the $\tilde{G}$-biharmonicity of a unit vector field (\ref{buvcond}) for the vector fields of the components of the frame field becomes:
\begin{center}
\begin{tabular}{c}
$T(e_{1})=\Big((\frac{a(a+d)}{\alpha}-\frac{b}{\varphi})(\frac{d}{\varphi}-4)-\frac{a+d}{\alpha}(\frac{\alpha +ad}{\varphi}-2a)\Big)be_{3}+2d\Big(\frac{a(a+d)}{\alpha}-\frac{b}{\varphi}\Big)e_{1},$\\\\
$T(e_{2})=\frac{ab}{\alpha \varphi}(d^{2}+\alpha-ad)e_{3}+\Big[\frac{b^{2}}{\varphi}+\frac{ad}{\alpha}(a+d)\Big]e_{2}$,\\\\
$T(e_{3})=\frac{8}{\varphi}b^{2}\Big(2-\frac{ad}{\alpha}\Big)e_{3}$.
\end{tabular}
\end{center}

\subsubsection*{Biharmonicity of $e_1$:}

From the expression of $T(e_{1})$ above,  $e_{1}$ is $\tilde{G}$-biharmonic if and only if
\begin{equation}\label{bih-e1}
  \left\{\begin{array}{l}
           d\Big(\frac{a(a+d)}{\alpha}-\frac{b}{\varphi}\Big)=0,  \\
           b\Big((\frac{a(a+d)}{\alpha}-\frac{b}{\varphi})(\frac{d}{\varphi}-4)-\frac{a+d}{\alpha}(\frac{\alpha +ad}{\varphi}-2a)\Big)=0
         \end{array}
  \right.
\end{equation}
According to the first equation of \eqref{bih-e1}, we have two cases:

\textbf{Case 1:} $d=0$. In this case, the second equation of \eqref{bih-e1} is equivalent to either $b=0$ or $b=\frac{a}{4\alpha}(2a(a+c) +\alpha)$.

\textbf{Case 2:} $d \neq 0$. In this case, the first equation of \eqref{bih-e1} is equivalent to
\begin{equation}\label{bih-e1-1}
  b=\frac{a(a+c)\varphi}{\alpha},
\end{equation}
and the second equation is equivalent to
\begin{equation}\label{bih-e1-2}
  \frac{a+d}{\alpha}\Big(\frac{\alpha +ad}{\varphi}-2a\Big)b=0.
\end{equation}
Taking into account \eqref{bih-e1-1}, equation \eqref{bih-e1-2} is equivalent to
\begin{itemize}
  \item either $a+d=b=0$. In this case, $e_1$ is a harmonic map $(M,g) \rightarrow (T_1M,\tilde{G})$;
  \item or $b \neq 0$, $a+d \neq 0$ and $\frac{\alpha +ad}{\varphi}-2a=0$. Then a simple calculation gives
  \begin{equation}\label{bih-e1-3}
    b^2=-a\varphi=\frac{-\alpha -ad}{2}.
  \end{equation}
  In particular, we have $a \neq 0$. Substituting from \eqref{bih-e1-3} into \eqref{bih-e1-1}, we obtain by virtue of $b \neq 0$
  \begin{equation}\label{bih-e1-4}
    b=-\frac{\alpha}{a+d}.
  \end{equation}
  Comparing \eqref{bih-e1-3} and \eqref{bih-e1-4}, we obtain $2\alpha^2 +(a+d)^2\alpha +a(a+d)^2d=0$. Noticing that this quadratic equation in $\alpha$ has a non-zero solution if and only if $d \in ]-\infty,3a-2\sqrt{2}|a|] \cup [3a+2\sqrt{2}|a|,+\infty[$ with $d \neq 0$. In this case the solutions are given by
  \begin{equation}\label{bih-e1-5}
    \alpha=\frac{-|a+d|}{4} (|a+d| \pm \sqrt{(a+d)^2 -8ad}).
  \end{equation}
\end{itemize}

\begin{Rq}
   Equation \eqref{bih-e1-5} shows that $\alpha$ depends on 2 parameters $a$ and $d$. On the other hand, from \eqref{bih-e1-4} we deduce that $b$ depends on $a$ and $d$. Consequently, by $a(a+c) -b^2=\alpha$, it is clear that $c$ depends on $a$ and $d$.
\end{Rq}

\subsubsection*{Biharmonicity of $e_2$:}

From the expression of $T(e_{2})$ above,  $e_{2}$ is $\tilde{G}$-biharmonic if and only if
\begin{equation}\label{bih-e2}
  \left\{\begin{array}{l}
           \frac{ab}{\alpha \varphi}(d^{2}+\alpha-ad)=0,  \\
           \frac{b^{2}}{\varphi}+\frac{ad}{\alpha}(a+d)=0.
         \end{array}
  \right.
\end{equation}
According to the first equation of \eqref{bih-e2}, we have two cases:

\textbf{Case 1:} $b=0$. In this case, since $a \neq 0$ ($\alpha \neq 0$), the second equation of \eqref{bih-e2} is equivalent to either $d=0$ or $a+d=0$.

\textbf{Case 2:} $b \neq 0$. In this case, the first equation of \eqref{bih-e2} is equivalent to either $a =0$ or $\alpha=d(a-d)$. But $a=0$ gives by the second equation of \eqref{bih-e2} $b=0$, and then $\alpha=0$, which is a contradiction. We deduce that
\begin{equation}\label{bih-e2-1}
  \alpha=d(a-d)
\end{equation}
and, in particular, $d \neq 0$ and $d \neq a$ (since $\alpha \neq 0$). Using \eqref{bih-e2-1}, the second equation of \eqref{bih-e2} is equivalent to
\begin{equation}\label{bih-e2-2}
  b^2=-\frac12 (a+d)(2a-d).
\end{equation}
In particular, $(a+d)(2a-d)<0$.

\subsubsection*{Biharmonicity of $e_3$:}

From the expression of $T(e_{3})$ above,  $e_{3}$ is $\tilde{G}$-biharmonic if and only if
\begin{itemize}
  \item either $b=0$. In this case, $e_{3}$ is a harmonic map $(M,g) \rightarrow (T_1M,\tilde{G})$;
  \item or $b \neq 0$ and $2\alpha=ad$. In particular, $a\neq 0$ and $d \neq 0$.
\end{itemize}
To summarize the previous discussion, we can state the following

\begin{Prop}
  Let \textit{Sol}$_{3}$ be the solvable Lie group ${\mathbb{R}}^{3}$ equipped with the metric  $g_{sol}=e^{2z}dx^{2}+e^{-2z}dy^{2}+dz^2$ and $\tilde{G}$ be an arbitrary pseudo-Riemannian $g$-natural metric on its unit tangent bundle. Let $\{e_1,e_2,e_3\}$ be the orthonormal frame field given by
\begin{equation*}
e_{1}=e^{-z}\frac{\partial}{\partial x}, \phantom {XXX}e_{2}=e^{z}\frac{\partial}{\partial y},\phantom {XXX}e_{3}=\frac{\partial}{\partial z}.
\end{equation*}
\begin{enumerate}
  \item $e_1$ is a $\tilde{G}$-biharmonic unit vector field if and only if one of the following assertions holds
  \begin{itemize}
    \item $d=b=0$ and $a \neq 0$. In this case, $e_1$ is proper $\tilde{G}$-biharmonic.
    \item $d=0$, $b \neq 0$ and $b=\frac{a}{4\alpha}(2a(a+c) +\alpha)$. In this case, $e_1$ is proper $\tilde{G}$-biharmonic.
    \item $d \neq 0$ and $b=a+d=0$. In this case, $e_1$ is a harmonic map $(M,g) \rightarrow (T_1M,\tilde{G})$.
    \item $a \neq 0$, $d \in \Big(]-\infty,3a-2\sqrt{2}|a|] \cup [3a+2\sqrt{2}|a|,+\infty[\Big) \setminus \{0\}$ and $b$ and $\alpha$ are given by \eqref{bih-e1-4} and \eqref{bih-e1-5}, respectively. In this case, $e_1$ is proper $\tilde{G}$-biharmonic.
  \end{itemize}
  \item $e_2$ is a $\tilde{G}$-biharmonic unit vector field if and only if one of the following assertions holds
  \begin{itemize}
    \item $d=b=0$ and $a \neq 0$. In this case, $e_2$ is proper $\tilde{G}$-biharmonic.
    \item $b=0$, $a+d=0$. In this case, $e_2$ is a harmonic map $(M,g) \rightarrow (T_1M,\tilde{G})$.
    \item $a \neq 0$, $d \neq 0$, $(a+d)(2a-d)<0$ and $\alpha$ and $b$ are given by \eqref{bih-e2-1} and \eqref{bih-e2-2}, respectively. In this case, $e_2$ is proper $\tilde{G}$-biharmonic.
  \end{itemize}
  \item $e_3$ is a $\tilde{G}$-biharmonic unit vector field if and only if one of the following assertions holds
  \begin{itemize}
    \item $b=0$. In this case, $e_3$ is a harmonic map $(M,g) \rightarrow (T_1M,\tilde{G})$.
    \item $a \neq 0$, $d \neq 0$, $b \neq 0$, $a+d=0$ and $2\alpha=ad$. In this case, $e_2$ is proper $\tilde{G}$-biharmonic.
  \end{itemize}
\end{enumerate}
\end{Prop}

\end{Ex}

\begin{Ex}\label{Ex-Sp-unitary}

On $SU(2)$ equipped with a left-invariant Riemannian metric, take an orthonormal basis $\{e_{1},e_{2},e_{3}\}$ of the Lie algebra $su(2)$ such that\cite{su2}:
$$[e_{1},e_{2}]=\lambda_{3} e_{3},\:\:\: [e_{2},e_{3}]=\lambda_{1} e_{1},\:\:\:[e_{3},e_{1}]=\lambda_{2} e_{2},$$

where $\lambda_{1}, \lambda_{2}$ and $\lambda_{3}$ are strictly positive constants with $\lambda_{1} \geq \lambda _{2} \geq \lambda_{3}$.

Then the Levi-Civita connexion $\nabla$ is determined by:
\begin{center}
\begin{tabular} {ccc}
$\nabla_{e_{1}}e_{1}=0 $,&$\nabla_{e_{1}}e_{2}=\mu_{1} e_{3}$,&$\nabla_{e_{1}}e_{3}=-\mu_{1} e_{2}$,\\\\
$\nabla_{e_{2}}e_{1}=-\mu_{2} e_{3} $,&$\nabla_{e_{2}}e_{2}=0$,&$\nabla_{e_{2}}e_{3}=\mu_{2} e_{1}$,\\\\
$\nabla_{e_{3}}e_{1}=\mu_{3} e_{2} $,&$\nabla_{e_{3}}e_{2}=-\mu_{3} e_{1}$,&$\nabla_{e_{3}}e_{3}=0$.\\\\
\end{tabular}
\end{center}
where: $$\mu_{i}=\frac{1}{2}(\lambda_{1}+\lambda_{2}+\lambda_{3})-\lambda_{i}, \:\:i=1,2,3.$$
A further computation gives:
\begin{center}
\begin{tabular}{ccc}
$R(e_{1},e_{2})e_{2}=(\lambda_{3}\mu_{3}-\mu_{1}\mu_{2})e_{1}$,&$R(e_{1},e_{3})e_{3}=(\lambda_{2}\mu_{2}-\mu_{1}\mu_{3})e_{2}$,&$R(e_{2},e_{1})e_{1} =(\lambda_{3}\mu_{3}-\mu_{1}\mu_{2})e_{2}$,\\\\
$R(e_{2},e_{3})e_{3}=(\lambda_{1}\mu_{1}-\mu_{2}\mu_{3})e_{2}$,&$R(e_{3},e_{1})e_{1}=(\lambda_{2}\mu_{2}-\mu_{1}\mu_{3})e_{3}$,&$R(e_{3},e_{2})e_{2} =(\lambda_{1}\mu_{1}-\mu_{2}\mu_{3})e_{3}$,\\\\
\end{tabular}
\end{center}
and the other components are zero. Furthermore, we get:

$$\underset{i=1}{\overset{m}{\sum}}\Big\{2R(e_{i},e_{j})\nabla_{e_{i}}e_{j}+(\nabla_{e_{i}}R)(e_{i},e_{j})e_{j}+R(e_{i},\nabla_{e_{i}}e_{j})e_{j}\Big\}=0,\:\:\: j=1,2,3.$$

$$div(e_{i})=S(e_{i})=\tilde{\tau}_{v}(e_{i})=0, \:\:\: i=1,2,3.$$
\begin{center}
\begin{tabular}{lcr}
$Q(e_{1})=(\lambda_{3}\mu_{3}-\mu_{1}\mu_{2}+\lambda_{2}\mu_{2}-\mu_{1}\mu_{3})e_{1}$,& &$\overline{\Delta}e_{1}=(\mu_{2}^{2}+\mu_{3}^{2})e_{1}$, \\\\
$Q(e_{2})=(\lambda_{3}\mu_{3}-\mu_{1}\mu_{2}+\lambda_{1}\mu_{1}-\mu_{2}\mu_{3})e_{2}$,& &$\overline{\Delta}e_{2}=(\mu_{1}^{2}+\mu_{3}^{2})e_{2}$,\\\\
$Q(e_{3})=(\lambda_{2}\mu_{2}-\mu_{1}\mu_{3}+\lambda_{1}\mu_{1}-\mu_{2}\mu_{3})e_{3}$,& &$\overline{\Delta}e_{3}=(\mu_{1}^{2}+\mu_{2}^{2})e_{3}$.
\end{tabular}
\end{center}


To simplify the case we can suppose that $\lambda_{1}=\lambda_{2}=\lambda$ then $\mu_{1}=\frac{\lambda_{3}}{2}$, $\mu_{2}=\frac{\lambda_{3}}{2}$, $\mu_{3}=\lambda-\frac{\lambda_{3}}{2}$ and we obtain
\begin{center}
\begin{tabular}{ccc}
$\tilde{\tau}_{h}(e_{1})=-\frac{b}{\varphi}(\lambda-\lambda_{3})^{2}e_{1}$,&$\tilde{\tau}_{h}(e_{2})=-\frac{b}{\varphi}(\lambda-\lambda_{3})^{2}e_{2}$,&$\tilde{\tau}_{h}(e_{3}) =-\frac{b}{\varphi}\frac{\lambda_{3}^{2}}{2} e_{3}$,
\end{tabular}
\end{center}
Using the formulas above, the condition of the biharmonicity of a unit vector field (\ref{buvcond}) for the vector fields of the frame field becomes:
\begin{center}
\begin{tabular}{c}
$T(e_{1})=-\frac{b^{2}}{\varphi}(\lambda-\lambda_{3})^{2}\Big[\lambda^{2}-\lambda\lambda_{3}(1-\frac{ad}{\alpha})-\lambda_{3}^{2}(-1+\frac{ad}{2\alpha}) \Big]e_{1}$,\\\\
$T(e_{2})=-\frac{b^{2}}{\varphi}(\lambda-\lambda_{3})^{2}\Big[\lambda^{2}-\lambda\lambda_{3}(1-\frac{ad}{\alpha})-\lambda_{3}^{2}(-1+\frac{ad}{2\alpha}) \Big]e_{2}$,\\\\
$T(e_{3})= 0$.
\end{tabular}
\end{center}
So $e_{3}$ is a harmonic map $(M,g) \rightarrow (T_1M,\tilde{G})$.\\\\
But $e_{1}$ and $e_{2}$ are $\tilde{G}$-biharmonic unit vector fields if and only if
\begin{itemize}
\item either $b=0$. In this case, $e_{1}$ and $e_{2}$ are harmonic maps $(M,g) \rightarrow (T_1M,\tilde{G})$;
\item or $\lambda_{3}=\lambda$. In this case, $e_{1}$ and $e_{2}$ are also harmonic maps $(M,g) \rightarrow (T_1M,\tilde{G})$;
\item or $\lambda^{2}-\lambda\lambda_{3}(1-\frac{ad}{\alpha})-\lambda_{3}^{2}(-1+\frac{ad}{2\alpha})=0$. Put $\lambda=\sigma \lambda_3$, so that $\sigma \geq 1$, and the quadratic equation above becomes $\sigma^{2}-(1-\frac{ad}{\alpha})\sigma-(-1+\frac{ad}{2\alpha})=0$, whose solutions are given by
    \begin{equation}\label{bih-SU-1}
      \sigma= \frac{1}{2\alpha}\Big(\alpha -ad \pm \sqrt{a^2d^2 -3\alpha^2}\Big),
    \end{equation}
    under the condition $a^2d^2 \geq 3\alpha^2$.

    \textbf{Case 1:} $\sigma= \frac{1}{2\alpha}\Big(\alpha -ad + \sqrt{a^2d^2 -3\alpha^2}\Big)$. This yields $\alpha(1-2\sigma)-ad= \sqrt{a^2d^2 -3\alpha^2}$, under the condition $\alpha(1-2\sigma)-ad \geq 0$. A routine calculation gives
    \begin{equation}\label{bih-SU-2}
      ad= \frac{2(\sigma^2 -\sigma +1)}{1-2\sigma}\alpha.
    \end{equation}
    Note that the condition $\sigma \geq 1$ in equation \eqref{bih-SU-2} gives the condition $a^2d^2 \geq 3\alpha^2$.

    \textbf{Case 2:} $\sigma= \frac{1}{2\alpha}\Big(\alpha -ad - \sqrt{a^2d^2 -3\alpha^2}\Big)$. The same arguments as in case 1, gives \eqref{bih-SU-2}, with the condition $\alpha(1-2\sigma)-ad \leq 0$.
\end{itemize}

Summarizing the last discussion, we get
\begin{Prop}
  Let $SU(2)$ equipped with a left-invariant Riemannian metric and $\tilde{G}$ be an arbitrary pseudo-Riemannian $g$-natural metric on its unit tangent bundle. Let $\{e_{1},e_{2},e_{3}\}$ be an orthonormal basis of the Lie algebra $su(2)$ such that
$$[e_{1},e_{2}]=\sigma \lambda e_{3},\:\:\: [e_{2},e_{3}]=\lambda e_{1},\:\:\:[e_{3},e_{1}]=\lambda e_{2},$$
where $\lambda>0$ and $\sigma\geq 1$. Then
\begin{enumerate}
  \item $e_{3}$ is a harmonic map $(M,g) \rightarrow (T_1M,\tilde{G})$.
  \item $e_1$ (resp. $e_2$) is a $\tilde{G}$-biharmonic unit vector field if and only if one of the following assertions holds
  \begin{itemize}
    \item $b=0$ or $\sigma=1$. In this case, $e_1$ and $e_2$ are harmonic maps $(M,g) \rightarrow (T_1M,\tilde{G})$.
    \item $b \neq 0$, $\sigma \neq 1$ and equation \eqref{bih-SU-2} holds. In this case $e_1$ and $e_2$ are proper $\tilde{G}$-biharmonic unit vector fields.
  \end{itemize}
\end{enumerate}
\end{Prop}
\end{Ex}

\subsection{Some particular cases}

    Next we will investigate some particular cases of Theorem \ref{buvcond}:

    \begin{The}
    Let $(M,g)$ be a Riemannian manifold and $T_{1}M$ be its unit tangent bundle equipped with an arbitrary pseudo-Riemannian Kaluza-Klein metric $\tilde{G}$ (i.e. $b=d=0$). Let $U \in \mathfrak{X}_{1}(M)$.  $U$ is a $\tilde{G}$-biharmonic unit vector field if and only if
     \begin{equation}
     \begin{array}{rcl}
     &&\frac{a}{a+c}\underset{i=1}{\overset{m}{\sum}}\Big\{2R(e_{i},S(U))\nabla_{e_{i}}U+(\nabla_{e_{i}}R)(e_{i},S(U))U+R(e_{i},\nabla_{e_{i}}S(U))U\Big\}\\\\
          &&+\overline{\Delta}\overline{\Delta}U-\overline{\Delta}\big[g(\overline{\Delta}U,U)U\big]-g(\overline{\Delta}U,U) \overline{\Delta}U+g(\overline{\Delta}U,U)^{2}U=0.\nonumber
     \end{array}
     \end{equation}
     If $a=1$ and $c=0$ i.e. $\tilde{G}=\widetilde{g^{s}}$, the result above is the same derived by M. Markellos and H.Urakawa in \cite{Ur.unit} in the case of Sasaki metric.
    \end{The}

    Given a compact Riemannian manifold $(M,g)$. In \cite{Ab.unit}, the authors defined a \emph{haramonic unit vector field} $U$ on $M$ as a critical point of the energy functional $E(U)=\frac{1}{2}\int_{M}||dU||^{2}v_{g}$ restricted to $\mathfrak{X}^{1}(M)$, the set of all unit vector fields, and they showed that this kind of harmonicity of $U$ does not depend upon the choice of the chosen $g$-natural metric $\tilde{G}$. More precisely, $U$ is a harmonic unit vector field if and only if $\overline{\Delta}U$ and $U$ are collinear. As a consequence of this result we have

    \begin{Cr}
    Let $(M,g)$ be a Riemannian manifold and $T_{1}M$ be its unit tangent bundle equipped with an arbitrary pseudo-Riemannian Kaluza-Klein metric $\tilde{G}$ (i.e. $b=d=0$). A harmonic unit vector field $U$ on $M$ is a $\tilde{G}$-biharmonic unit vector field if and only if
         \begin{equation}
    \frac{a^{2}}{\alpha}\underset{i=1}{\overset{m}{\sum}}\Big\{2R(e_{i},S(U))\nabla_{e_{i}}U+(\nabla_{e_{i}}R)(e_{i},S(U))U+R(e_{i},\nabla_{e_{i}}S(U))U\Big\} \:\: and \:\: U \:\: are \:\:  collinear.
         \end{equation}
          If $a=1$ and $c=0$ i.e. $\tilde{G}=\widetilde{g^{s}}$, the result above is the same derived by M. Markellos and H. Urakawa in \cite{Ur.unit} in the case of Sasaki metric.
    \end{Cr}

    \begin{The}
    Let $(M,g)$ be a Riemannian manifold and $T_{1}M$ be its unit tangent bundle equipped with an arbitrary pseudo-Riemannian Kaluza-Klein type metric $\tilde{G}$ (i.e. $b=0$). Let $U \in \mathfrak{X}^{1}(M)$ satisfy $S(U)=-\frac{d}{a}\textup{div}(U)U$. Then $U$ is a $\tilde{G}$-biharmonic unit vector field if and only if
      \begin{equation}
      \begin{array}{rcl}
      &&\underset{i=1}{\overset{m}{\sum}}\Big\{2R(e_{i},\tilde{\tau}_{h}(U))\nabla_{e_{i}}U +(\nabla_{e_{i}}R)(e_{i},\tilde{\tau}_{h}(U))U+R(e_{i},\nabla_{e_{i}}\tilde{\tau}_{h}(U))U)+\frac{d}{a}g\big(\nabla_{e_{i}}U,\tilde{\tau}_{h}(U)\big)e_{i}\Big\}\\\\
         &&+\overline{\Delta}\tilde{\tau}_{v}(U) -\frac{d}{a}\nabla_{U}\tilde{\tau}_{h}(U)-g(\overline{\Delta}U,U)\tilde{\tau}_{v}(U)-\frac{d(a+c)}{a\varphi}\textup{div}(U)\tilde{\tau}_{h}(U)=0.
      \end{array}
      \end{equation}
    \end{The}

    \begin{Cr}
    Let $(M,g)$ be a Riemannian manifold of constant sectional curvature $k$ and $T_{1}M$ be equipped with an arbitrary pseudo-Riemannian Kaluza-Klein type metric $\tilde{G}$ (i.e. $b=0$) such that $d=-ak$ . A unit geodesic vector field $U$, i.e. $\nabla_{U}U=0$, is a $\tilde{G}$-biharmonic unit vector field if and only if
     \begin{equation}
      \begin{array}{rcl}
      &&2k\nabla_{U}\tilde{\tau}_{h}(U)+2k\:\textup{div}(U)\tilde{\tau}_{h}(U)+
      \frac{k(a+c)}{\varphi}\textup{div}(U)\tilde{\tau}_{h}(U)-\overline{\Delta}\tilde{\tau}_{v}(U)+g(\overline{\Delta}U,U)\tilde{\tau}_{v}(U)=0.
      \end{array}
      \end{equation}
    \end{Cr}

We can also discuss Theorem \ref{buvcond} for some special classes of vector fields such that \textbf{Killing vector fields} and \textbf{Reeb vector fields}.

Let $(M, g)$ be a Riemannian manifold. As it is well known, $U \in \mathfrak{X}(M)$ is a Killing vector field if the local 1-parameter group of $U$ consists of local isometries of $g$. Moreover, a vector field $U$ is a Killing vector field if and only if $\mathcal{L}_{U} g = 0$, where $\mathcal{L}$ denotes the Lie derivative, and we have in this case \cite{Ab.unit}:\\
$$div(U)=0,\:\:\: \nabla_{U}U=0,\:\:\:\overline{\Delta}U=QU,\:\:\:g(QU,U)=||\nabla U||^{2},$$
Let $\tilde{G}$ be an arbitrary pseudo-Riemanniann $g$-natural metric on $T_{1}M$ with $b=0$ and suppose that $S(U)=0$. Then the horizontal and vertical parts of the tension field are given by:
$$\tilde{\tau}_{h}(U)=0,\:\:\:\tilde{\tau}_{v}(U)=QU-||\nabla U||^{2}U=\mathcal{A}U. $$

Then, as a corollary of Theorem \ref{buvcond}, we get
    \begin{Cr}
    Let $(M,g)$ be a Riemannian manifold and $T_{1}M$ be equipped with an arbitrary pseudo-Riemannian $g$-natural Kaluza-Klein type metric $\tilde{G}$ (i.e. $b=0$). Let $U$ be a unit Killing vector field on $M$ which satisfies $S(U)=0$ and suppose that $QU$ and $U$ are not collinear. Then $U$ is a proper $\tilde{G}$-biharmonic unit vector field if and only if

      \begin{equation}
       \overline{\Delta}\mathcal{A}U=||\nabla U||^{2}   \mathcal{A}U.
       \end{equation}
    \end{Cr}

 Now we consider the special case when the unit vector filed is the Reeb vector field $\xi$ of a contact metric manifold $(M,\eta,g)$.\\

A $(2m+1)$-dimensional manifold $M$ is said to be a \textit{contact metric manifold} if it admits a global $1$-form $\eta$ such that $\eta \wedge (d\eta)^{m}\neq 0$. There is a unique vector field $\xi$, called the \textit{Reeb vector field} (or the characteristic vector field), such that $\eta(\xi)=1$ and $d\eta(\xi,.)=0$. Furthermore, a Riemannian metric $g$ is said to be an \textit{associated metric} if there exists a tensor field $\varphi$ of type $(1,1)$ such that
\begin{equation}\label{contact}
\eta=g(\xi,.),\:\:d\eta=g(.,\varphi .),\:\: \varphi^{2}=-I+\eta\otimes \xi.
\end{equation}
By (\ref{contact}) it follows that $\xi$ is a unit vector field on $(M,g)$, and it satisfies

\begin{center}
\begin{tabular}{c}
$\nabla \xi=-\varphi-\varphi h$, $\:\:\:\nabla_{\xi}\xi=0$, $\:\:\:div(\xi)=0$\\\\
$||\nabla \xi||^{2}=2m+tr(h^{2})=4m-g(Q\xi,\xi)$,$\:\:\:\overline{\Delta}\xi=4m\xi-Q\xi $
\end{tabular}
\end{center}
where $h=\frac{1}{2}\mathcal{L}_{\xi}\varphi$.

Let $\tilde{G}$ be an arbitrary pseudo-Riemannian $g$-natural metric on $T_{1}M$ with $b=0$, and we suppose that $S(\xi)=0$, then the horizontal and vertical parts of the tension filed are given by:
$$\tilde{\tau}_{h}(\xi)=0,\:\:\:\tilde{\tau}_{v}(\xi)=Q\xi-g(Q\xi,\xi)\xi=\mathcal{A}\xi. $$

Then we have

     \begin{Cr}
    Let $(M,\eta,g)$ be a contact metric manifold and $T_{1}M$ be equipped with an arbitrary pseudo-Riemannian $g$-natural Kaluza-Klein type metric $\tilde{G}$ (i.e. $b=0$). Suppose that the Reeb vector field $\xi$ on $M$ satisfies $S(\xi)=0$ and we suppose that $Q\xi$ and $\xi$ are not collinear. Then $\xi$ is a proper $\tilde{G}$-biharmonic unit vector field if and only if
       \begin{equation}
   \overline{\Delta}\mathcal{A}\xi=||\nabla \xi||^{2}   \mathcal{A}\xi.
        \end{equation}
     \end{Cr}

\newpage
\appendix
\makeatletter
\newcommand{\section@cntformat}{Appendix \thesection:\ }
\makeatother


\section{$F$-tensor fields}


In this Appendix, we give basic facts on $F$-tensor fields which will be needed in this work. For the proofs and more details, we refer to \cite{Ab.heri}.

\emph{$F$-tensor fields} are mappings $A :TM \oplus \underbrace{TM \oplus ...\oplus TM}_{s \; times} \rightarrow \bigsqcup_{x \in M} \otimes^r M_x$ which are linear in the last $s$ summands such that $\pi_2 \circ A =\pi_1 $, where $\pi_1$ and $\pi_2$ are the natural projections of the source and target fiber bundles of $A$ respectively. For $r=0$ and
$s=2$, we obtain the classical notion of \emph{$F$-metrics}. So, if we denote by $\oplus$ the fibered product of fibered manifolds, then $F$-metrics are mappings $TM \oplus TM \oplus TM \to \mathbb{R}$ which are linear in the second and the third argument.\par

Fix $(x,u) \in TM$ and a system of normal coordinates $S:=(U;x^1,..., x^m)$ of $(M,g)$ centered at $x$. Then we can define on $U$ the vector field $\mathrm{U}:= \sum_i u^i \;
\frac{\partial}{\partial x^i}$, where $(u^1,...,u^m)$ are the coordinates of $(x,u)$ with respect to the basis $((\frac{\partial}{\partial x^i})_x;\; i=1,...,m)$ of $M_x$. \par

Let $P$ be an $F$-tensor field of type $(p,q)$ on $M$. Then, on $U$, we can define a $(p,q)$-tensor field $P_u^S$ (or $P_u$ if there is no risk of confusion), associated to $u$ and $S$, by
\begin{equation}\label{tens-ass}
P_u(X_1,...,X_q):= P(\mathrm{U}_z;X_1,...,X_q),
\end{equation}
for all $(X_1,...,X_q) \in M_z$, $z \in U$. Informally, we can say that we have \textgravedbl\emph{tensorized}\textacutedbl $P$ at $u$ with respect to $S$.\par

On the other hand, if we fix $x \in M$ and $q$ vectors $X_1,...,X_q$ in $M_x$, then we can define a $C^\infty$-mapping $P_{(X_1,...,X_q)}:M_x \rightarrow \otimes^p M_x$ , associated to $(X_1,...,X_q)$, by
\begin{equation}\label{map-ass}
P_{(X_1,...,X_q)}(u) :=P(u;X_1,...,X_q),
\end{equation}
for all $u \in M_x$.\par

Let $s>t$ be two non negative integers, $T$ be a $(1,s)$-tensor field on $M$ and $P^T$ be an $F$-tensor field, of type $(1,t)$, of the form
\begin{equation}\label{tens}
P^T(u;X_1,...,X_t) =T(X_1,...,u,...,u,...,X_t),
\end{equation}
for all $(u,X_1,...,X_t) \in TM \oplus ... \oplus TM$, i.e. $u$ appears $s-t$ times at positions $i_1,...,i_{s-t}$ in the expression of $T$. Then

\begin{itemize}
  \item [*] $P^T_u$ is a $(1,t)$-tensor field on a neighborhood $U$ of $x$ in $M$, for all $u \in M_x$;
  \item [*] $P^T_{(X_1,...,X_t)}$ is a $C^\infty$-mapping $M_x  \rightarrow M_x$, for all $X_1,...,X_t$ in $M_x$.
\end{itemize}
Furthermore, we have

\begin{Lm}
1) The covariant derivative of $P^T_u$, with respect to the Levi-Civita connection of $(M,g)$, is given by:
\begin{equation}\label{der-ass}
(\nabla_{X} P^T_u)(X_1,...,X_t) = (\nabla _X T)(X_1,...,u,...,u ,..., X_{t}) ,
\end{equation}
for all vectors $X$, $X_1$,..., $X_{t}$ in $M_x$, where $u$ appears at positions $i_1,...,i_{s-t}$ in the right hand side of the preceding formula.                                        \\
2) The differential of $P^T_{(X_1,...,X_t)}$, at $u \in M_x$, is given by: \arraycolsep1.5pt
\begin{eqnarray}
d_u(P^T_{(X_1,...,X_t)})(X)
         & = & T(X_1,...,X,...,u,...,X_{t}) \label{diff-ass}\\
         &   & +... +T(X_1,...,u,...,X,...,X_{t}), \nonumber
\end{eqnarray}
\arraycolsep5pt for all $X \in M_x$.
\end{Lm}

We have also the following:

\begin{Lm}\label{der-T}Let $T$ be a $(1,s)$-tensor field on $M$. Then
\arraycolsep1.5pt
\begin{eqnarray*}
1)\bar \nabla_{X^h}h\{ T  ( X_1 ,...,u,...,u,...,X_{t})\}
         & = & h\{(\nabla_{X} P^T_u)((X_1)_x,...,(X_t)_x) \\
         &   & +A(u;X,T_x(X_1,...,u,...,u,...,X_{t}))\}\\
         &   & + v\{ B(u;X,T_x(X_1,...,u,...,u,...,X_{t}))\},
 \\
2) \bar \nabla_{X^v}h\{  T  (X_1,...,u,...,u,...,X_{t})\}
         & = & h\{d_u(P^T_{((X_1)_x,...,(X_t)_x)})(X)\\
         &   & +C(u;T_x(X_1,...,u,...,u,...,X_{t}),X)\}\\
         &   & +  v\{ D(u;T_x(X_1,...,u,...,u,...,X_{t}),X)\},
\\
3) \bar \nabla_{X^h}v\{  T  (X_1,...,u,...,u,...,X_{t})\}
       & = & h\{C(u;X,T_x(X_1,...,u,...,u  ,...,X_{t}))\} \\
       &   & +v\{(\nabla_{X} P^T_u)((X_1)_x,...,(X_t)_x) \\
       &   & +D(u;X,T_x(X_1,...,u,...,u  ,...,X_{t}))\},
  \\
4) \bar \nabla_{X^v}v\{ T  (X_1,...,u,...,u,...,X_{t})\}
  & = & h\{E(u;X,T_x(X_1,...,u,...,u,...,X_{t}))\} \\
  &   & + v\{d_u(P^T_{((X_1)_x,...,(X_t)_x)})(X)) \\
  &   & +  F(u;X,T_x(X_1,...,u,...,u,...,X_{t}))\},
\end{eqnarray*}
\arraycolsep5pt     for all vector fields $X_1$,..., $X_{t}$ on
$M$ and $X \in M_x$, where $u$ appears at positions
$i_1,...,i_{s-t}$ in in any expression of $T$. Here, $X^h$ and
$X^v$ are taken at $(x,u)$.
\end{Lm}

Now, let $P$ be the $F$-tensor field of type $(1,t)$ of the form
\begin{equation}\label{g-form-P}
  P(u;X_1,...,X_t) =\sum_i f_i^P (r^2) T_i(X_1,...,u,...,u,...,X_t),
\end{equation}
where $f_i: \mathbb{R}^+ \rightarrow \mathbb{R}$ are real-valued functions on $\mathbb{R}^+$, and any $T_i$ is a $(1,s_i)$-tensor field on $M$, $s_i > t$, with the $s_i$'s not necessarily equal. Then, we have

\begin{Lm}\label{der-diff1}
Let $P$ be an $F$-tensor field, of type $(1,t)$, on $M$ given by (\ref{g-form-P}). Then \arraycolsep0pt
\begin{eqnarray*}
1) (\nabla _X  P_u) (X_1,...,X_t) & = & \sum_i f_i^P (r^2) (\nabla_X T_i)(X_1,...,u,...,u,...,X_t),                      \\
2) d_u(P_{(X_1,...,X_t)}) (X) & = & 2 \sum_i (f_i^P)^\prime (r^2) g(X,u) T_i(X_1,...,u,...,u,...,X_t) \\
 & & +\sum_i f_i^P (r^2)\{ T_i(X_1,...,X,...,u,...,X_t)+... \\
 & & +T_i(X_1,...,u,...,X,...,X_t)\},
\end{eqnarray*}
for all $u,X,X_1,...,X_t \in M_x$.
\end{Lm}

If we denote by $h\{P(u;X_1,...,X_t)\}$ (resp. $v\{P(u;X_1,...,X_t)\}$) the quantity \arraycolsep0pt
\begin{eqnarray}
h\{P(u;X_1,...,X_t)\} & = & \sum_i f_i^P (r^2) h\{T_i(X_1,...,u,...,u,...,X_t)\}  \label{h-P}          \\
\mbox{(resp.}\;\; v\{P(u;X_1,...,X_t)\} & = & \sum_i f_i^P (r^2) v\{T_i(X_1,...,u,...,u,...,X_t)\}),  \label{v-P}
\end{eqnarray}
then we can assert
\begin{Lm}\label{der-F-T}
\arraycolsep1.5pt
\begin{eqnarray*}
1)\bar \nabla_{X^h}h\{P(u;X_1,...,X_t)\}
         & = & h\{(\nabla_{X} P_u)((X_1)_x,...,(X_t)_x) \\
         &   & +A(u;X,P(u;(X_1)_x,...,(X_t)_x))\}\\
         &   & + v\{ B(u;X,P(u;(X_1)_x,...,(X_t)_x))\},
 \\
2) \bar \nabla_{X^v} h\{P(u;X_1,...,X_t)\}
         & = & h\{d_u(P_{((X_1)_x,...,(X_t)_x)})(X)\\
         &   & +C(u;P(u;(X_1)_x,...,(X_t)_x),X)\}\\
         &   & +v\{ D(u;P(u;(X_1)_x,...,(X_t)_x),X)\},
\\
3) \bar \nabla_{X^h}v\{ P(u;X_1,...,X_t)\}
       & = & h\{C(u;X,P(u;(X_1)_x,...,(X_t)_x))\} \\
       &   & +v\{(\nabla_{X} P_u)((X_1)_x,...,(X_t)_x)\\
       &   & +D(u;X,P(u;(X_1)_x,...,(X_t)_x))\},
  \\
4) \bar \nabla_{X^v}v\{P(u;X_1,...,X_t)\}
       & = & h\{E(u;X,P(u;(X_1)_x,...,(X_t)_x))\} \\
       &   & + v\{d_u(P_{((X_1)_x,...,(X_t)_x)})(X))\\
       &   & +F(u;X,P(u;(X_1)_x,...,(X_t)_x))\},
\end{eqnarray*}
for all vector fields $X_1,...,X_t$ on $M$ and $X \in M_x$. Here $X^h$ and $X^v$ are taken at $(x,u)$.
\end{Lm}

Now, the existence of non-vanishing parallel vector field on a Riemannian manifold $(M,g)$  implies that all the $F$-tensors fields $A, B, C, D, E$ and $F$ reduce to the following form:
\begin{equation}\label{g-form-P2}
 P(u;X,Y)=f_{3}^{p}g(Y,u)X+f_{4}^{p}g(X,u)Y+f_{5}^{p}g(X,Y)u+f_{6}^{p}g(X,u)g(Y,u)u,
     \end{equation}
     wehre $f_{i}:\mathbb{R}^{+}\longrightarrow \mathbb{R}$ are real valued function on $\mathbb{R}^{+}$.
\begin{Lm}\label{der-diff2}
Let $P$ be an $F$-tensor field, of type $(1,2)$, on $M$ given by (\ref{g-form-P2}). Then \arraycolsep0pt
\begin{eqnarray*}
1) \nabla _X  P_u & = & 0,                      \\
2) d(P_{(X,Y)})_u (Z) & = & [f_{3}^{P}g(Y,Z)+2(f_{3}^{P})'g(Y,u)g(Z,u)]X \\
 & & +[f_{4}^{P}g(X,Z)+2(f_{4}^{P})'g(X,u)g(Z,u)]Y+f_{5}^{P}g(X,Y)Z\\
 & &+\{f_{6}^{P}[g(X,Z)g(Y,u)+g(Y,Z)g(X,u)]\\
  & &+2(f_{5}^{P})'g(X,Y)g(Z,u)+2(f_{6}^{P})'g(X,u)g(Y,u)g(Z,u)\} u,
\end{eqnarray*}
for all $u,X,Y,Z \in M_x$.
\end{Lm}


\end{document}